# LINEAR AND UNCONDITIONALLY ENERGY STABLE SCHEMES FOR THE BINARY FLUID-SURFACTANT PHASE FIELD MODEL

XIAOFENG YANG †⋆, LILI JU ‡

Abstract. In this paper, we consider the numerical solution of a binary fluid-surfactant phase field model, in which the free energy contains a nonlinear coupling entropy, a Ginzburg-Landau double well potential, and a logarithmic Flory-Huggins potential. The resulting system consists of two coupled, nonlinear Cahn-Hilliard type equations. We develop a set of first and second order time marching schemes for this system using the "Invariant Energy Quadratization" approach, in particular, the system is transformed into an equivalent one by introducing appropriate auxiliary variables and all nonlinear terms are then treated semi-explicitly. Both schemes are linear and lead to symmetric positive definite systems at each time step, thus they can be efficiently solved. We further prove that these schemes are unconditionally energy stable in the discrete sense. Various 2D and 3D numerical experiments are performed to validate the accuracy and energy stability of the proposed schemes.

## 1. Introduction

Surfactants are some organic compounds that can reduce the surface tension of the solution and allow for the mixing of immiscible liquids. A typical well-known example of immiscible liquids is the mixture of oil and water. There are many studies on the modeling and numerical simulations to investigate the binary fluid-surfactant system. In the pioneering work of Laradji et. al. [25, 26], the diffuse interface approach, or called the phase field method, was first used to study the phase transition behaviors of the monolayer microemulsion system, formed by surfactant molecules. The idea of the phase field method can be traced to Rayleigh [34] and Van der Waals [47] a century ago. Now this method has become a well-known effective modeling and simulation tool to resolve the motion of free interfaces between multiple material components. About the recent developments in advanced algorithms and computational technologies of the phase field method, we refer to [4–6, 9, 14, 22, 23, 27, 28, 30, 32, 44, 56] and references cited therein.

A variety of binary fluid surfactant phase field (BFS-PF) models had been well investigated in the past two decades. In [25, 26], two phase field variables are introduced to represent the local densities of the fluids, as well as the local concentration of the surfactant, respectively. There are two types of nonlinear energy terms in the model, including (i) the phenomenological Ginzburg-Landau (G-L) double well potential for the density variable to describe the phase separation behaviors of the fluid mixture, and (ii) the nonlinear coupling entropy term to ensure the high fraction of the surfactant near the fluid interface. Subsequently, the authors in [24] developed a modified model by adding an extra diffusion term and a G-L type potential for the concentration variable, in order to improve the stability. In [46] the logarithmic Flory-Huggins (F-H) potential was added in order to restrict the range of the concentration variable, while the nonlinear coupling entropy is essentially the same as that in [24–26]. A slightly different nonlinear coupling entropy was presented in [13], which could penalize the concentration to accumulate along the fluid interface. In [45], the authors further modified the model in [13] by adding the F-H potential for the local concentration variable as well.

⋆Corresponding author.
† Department of Mathematics, University of South Carolina, Columbia, SC 29208, USA. Email: xfyang@math.sc.edu. X. Yang's research is partially supported by the U.S. National Science Foundation under grant numbers DMS-1200487 and DMS-1418898.
‡Department of Mathematics, University of South Carolina, Columbia, SC 29208, USA. Email: ju@math.sc.edu. L. Ju's research is partially supported by the U.S. National Science Foundation under grant number DMS-1521965, and the U.S. Department of Energy under grant number DE-SC0008087-ER6539.

From the numerical point of view, it is very challenging to develop unconditionally energy stable schemes to discretize the nonlinear stiffness terms for the phase field type models, where the stiffness is originated from the thin interface thickness parameter. As a matter of fact, the simple fully-implicit or explicit type discretizations will induce very severe time step constraint (called conditionally energy stable) on the interfacial width [3,12,38], so they are not efficient in practice. Many efforts had been done in this direction in order to remove this type of time step constraint (cf. [6,11,12,18–20,29,31,33,35–38,40–43,49,49,50,52,53,58,60,61]). About these developed numerical techniques, we give a detailed discussion in Section 3. In addition, we emphasize that the "unconditional" here means the schemes have no constraints on the time step from stability point of view. However, large time step size will definitely induce large errors in practice. This fact motivates us to develop more accurate schemes, e.g., the second order time marching schemes while preserving the unconditional energy stability in this paper.

In addition to the stiffness issue from the interfacial width, we must notice that a particular specialty of the BFS-PF system is the strong *nonlinear couplings* between multiple phase field variables, that increases the complexity for algorithm developments to a large extent. Therefore, although a variety of phase field fluid surfactant models had been developed for over twenty years, there are very few successful attempts in designing efficient and energy stable schemes for them. Recently, in [15], Gu et. al. had developed a first order in time, nonlinear scheme to solve a particular BFS-PF model developed in [13,45] based on the convex splitting approach, where the convex part of the free energy potential is treated implicitly while the concave part is treated explicitly. Except an assumption that the approximate solutions always (accidentally) sit inside the domain of the logarithmic functional (such an assumption is often hard to hold in practical simulations), their arguments about the convex-concave decomposition for the coupling potential are not valid as well since it is not sufficient to justify the convexity of a function with multiple variables from the positivity of second order partial derivatives. In addition, their scheme is only first order in time, and its computational cost is relatively expensive due to the nonlinear nature.

Therefore, in this paper, the main purpose is to develop some more efficient and effective numerical schemes to solve the particular BFS-PF model that had been developed in [13, 45] since this model is a typical representative of nonlinear coupled multivariate fluid-surfactant models. We expect that our developed schemes can combine the following three desired properties, i.e., (i) *accurate* (ready for second order in time); (ii) *stable* (the unconditional energy dissipation law holds); and (iii) *easy to implement and efficient* (only need to solve some fully *linear* equations at each time step). To achieve such a goal, instead of using traditional discretization approaches like simple implicit, stabilized explicit, convex splitting, or other various tricky Taylor expansions to discretize the nonlinear potentials, we adopt the so-called *Invariant Energy Quadratization* (IEQ) method, which is a novel approach and had been successfully applied for other gradient flow models in the author's work (cf. [18,51,54,55,57,59]). But when this method is applied to the multivariate model such as the BFS-PF model considered in this paper, there are still some new challenges due to the nonlinear couplings between the multiple variables. Our essential idea is to transform the free energy into a quadratic form (since the nonlinear potential is usually bounded from below) of a set of new variables via a change of variables. The new, equivalent system still retains the similar energy dissipation law in terms of the new variables. One great advantage of such a reformulation is that all nonlinear terms can be treated semi-explicitly accordingly, which in turn produces a linear system. Moreover, the resulted linear operator of the system is symmetric positive definite, thus it can be solved by many available efficient linear solvers, e.g., CG, GMRES, or other Krylov subspace methods. We emphasize that the schemes developed in this paper is general enough to be extended to develop efficient and accurate schemes for a large class of gradient flow problems with multiple variables and/or complex nonlinearities in the free energy density.

The rest of the paper is organized as follows. In Section 2, the binary fluid-surfactant phase field model and its energy law are briefly introduced. In Section 3, we present our numerical schemes with respective first order and second order temporal accuracy for simulating the target model, and rigorously prove that the produced linear systems are symmetric positive definite and the schemes satisfy the unconditional energy stabilities. Various 2D and 3D numerical experiments are carried out in Section 4 to validate the accuracy and stability of these schemes. Finally, some concluding remarks are given in Section 5.

## 2. The BFS-PF model and its energy law

In the BFS system, monolayers of surfactant molecules form microemulsions as a random phase. Such microemulsion system usually exhibits various interesting microstructures, depending on the temperature or the composition. The phase field modeling approach simulates the dynamics of microphase separation in microemulsion systems using two phase field variables (order parameters). We now give the brief introduction of the BFS-PF model in [13, 45].

To label the local densities of the two fluids, fluid I and fluid II (such as water and oil), a phase field variable $\phi(\boldsymbol{x}, t)$ is introduced such that

$$\phi(\boldsymbol{x}, t) = \begin{cases} -1 & \text{fluid I}, \\ 1 & \text{fluid II}, \end{cases} \tag{2.1}$$

with a thin smooth transition layer of width $O(\epsilon)$ connecting the two fluids. Thus the interface of the mixture is described by the zero level set $\Gamma_t = \{\boldsymbol{x} : \phi(\boldsymbol{x}, t) = 0\}$. Let $F(\phi) = (\phi^2 - 1)^2$ be the classical G-L double well potential, one can define the mixing free energy functional as

$$\boldsymbol{E}_1(\phi) = \int_\Omega \Big(\frac{\epsilon}{2}|\nabla\phi|^2 + \frac{1}{4\epsilon}F(\phi)\Big)d\boldsymbol{x}, \tag{2.2}$$

where $\Omega$ denotes the computational domain of the problem. Note that $\boldsymbol{E}_1(\phi)$ is the most commonly used free energy in phase field models so far, where the first term in (2.2) contributes to the hydrophilic type (tendency of mixing) of interactions between the materials and the second term represents the hydrophobic type (tendency of separation) of interactions. As the consequence of the competition between the two types of interactions, the equilibrium configuration will include a diffusive interface. The systems induced from this part of energy have been well studied in [1, 17, 21, 28].

To represent the local concentration of surfactants, we take another phase field variable $\rho(\boldsymbol{x}, t)$ and assume its associated free energy as

$$\boldsymbol{E}_2(\rho) = \beta \int_\Omega G(\rho) d\boldsymbol{x}, \tag{2.3}$$

where $G(\rho) = \rho\ln\rho + (1-\rho)\ln(1-\rho)$ and $\beta$ is a positive parameter. We note that $G(\rho)$ is the logarithmic Flory-Huggins type energy potential which restricts the value of $\rho$ to be inside the domain of $(0, 1)$, and $\rho$ will reach its upper bound if the interface is fully saturated with surfactant (cf. [46]).

Lastly, due to the special property of surfactants that can alter the interfacial tension, the fraction of the surfactant staying at the fluid interface is high. Thus a local, nonlinear coupling entropy term between $\phi$ and $\rho$ is imposed as

$$\boldsymbol{E}_3(\phi, \rho) = \frac{\alpha}{2} \int_\Omega (\rho - |\nabla\phi|)^2 d\boldsymbol{x}, \tag{2.4}$$

where $\alpha$ is a positive parameter. This part of energy is the penalty term that enables the concentration to accumulate near the interface with a relatively high value.

The total energy of the system is then a sum of these three terms

$$\boldsymbol{E}_{tot}(\phi, \rho) = \boldsymbol{E}_1(\phi) + \boldsymbol{E}_2(\rho) + \boldsymbol{E}_3(\phi, \rho). \tag{2.5}$$

Using the variational approach, a nonlinear, coupled Cahn–Hilliard type phase field system for $(\phi, \rho)$ can be derived as

$$\phi_t = M_1 \Delta\mu_\phi, \tag{2.6}$$

$$\mu_\phi = -\epsilon\Delta\phi + \frac{1}{\epsilon}\phi(\phi^2 - 1) + \alpha\nabla\cdot\Big((\rho - |\nabla\phi|)\boldsymbol{Z}\Big), \tag{2.7}$$

$$\rho_t = M_2 \Delta\mu_\rho, \tag{2.8}$$

$$\mu_\rho = \alpha(\rho - |\nabla\phi|) + \beta\ln\Big(\frac{\rho}{1-\rho}\Big), \tag{2.9}$$

where $\boldsymbol{Z} = \frac{\nabla\phi}{|\nabla\phi|}$, $M_1$ and $M_2$ are the mobility parameters. For simplicity, we take the periodic boundary conditions to remove all complexities from the boundary integrals.

It is then straightforward to obtain the PDE energy law for BFS-PF Cahn-Hilliard system (2.6)-(2.9). Denote by $(f(\boldsymbol{x}), g(\boldsymbol{x})) = \int_\Omega f(\boldsymbol{x})g(\boldsymbol{x})d\boldsymbol{x}$ the $L^2$ inner product of any two functions $f(\boldsymbol{x})$ and $g(\boldsymbol{x})$, and by $\|f\| = \sqrt{(f,f)}$ the $L^2$ norm of the function $f(\boldsymbol{x})$. Let $L^2_{per}(\Omega)$ denote the subspace of all functions with the periodic boundary condition in $L^2(\Omega)$. By taking the sum of the $L^2$ inner products of (2.6) with $\mu_\phi$, of (2.7) with $\phi_t$, of (2.8) with $\mu_\rho$, and of (2.9) with $\rho_t$, we obtain

$$\frac{d}{dt}\boldsymbol{E}_{tot} = -M_1\|\nabla\mu_\phi\|^2 - M_2\|\nabla\mu_\rho\|^2 \le 0. \tag{2.10}$$

In the sequel, our goal is to design temporal approximation schemes which satisfy the discrete version of the continuous energy law (2.10).

**Remark 2.1.** *It is natural to incorporate hydrodynamic effects into the above model by introducing the extra stress induced from the free energy* (2.5) *using the similar approach as in [23, 28]. However, the hydrodynamics coupled system will present further numerical challenges, e.g., how to decouple the computations of the velocity from the phase variables. Since this paper is focused on the development of efficient linear schemes for solving the nonlinearly coupled Cahn-Hilliard equations with multiple variables, the details of numerical schemes for the hydrodynamics coupled model that are in the similar vein as [28, 29, 39, 40, 52], will be implemented in our future work.*

## 3. Numerical schemes

We next construct unconditionally energy stable, linear numerical schemes for solving the BFS-PF model (2.6)-(2.9). To this end, there are mainly three challenging issues, including (i) how to discretize the cubic term associated with the G-L double well potential; (ii) how to discretize the the logarithmic term induced by the F-H potential; and (iii) how to discretize the local coupling entropy terms associated with $\rho$ and $\phi$.

The discretization of the nonlinear, cubic polynomial term $f(\phi) = F'(\phi)$ had been well studied in many works (cf. [23, 28, 30, 32, 38, 44]). In summary, there are mainly two commonly used techniques to discretize $f$ in order to preserve the unconditional energy stability. The first one is the so-called convex splitting approach [11, 48], where the convex part of the potential is treated implicitly and the concave part is treated explicitly. The convex splitting approach is energy stable, however, it produces nonlinear schemes at most cases, thus the implementations are often complicated and the computational costs are high. Moreover, there are very few studies in the literature to apply the convex splitting approach to the case of multiple variables. The second technique is the so-called linear stabilization approach [6, 38, 49], where the nonlinear term is simply treated explicitly. In order to preserve the energy law, a linear stabilizing term has to be added, and the magnitude of that term usually depends on the upper bound of the second order derivative of the G-L potential. The stabilized approach introduces purely linear schemes, thus it is easy to implement and solve. However, the derivative usually does not have finite upper bound. A feasible remedy is to make some reasonable revisions to the nonlinear potential in order to obtain a finite bound, for example, the quadratic order cut-off functions for the G-L potential (cf. [6, 38, 49]). Such method is particularly reliable for those models with maximum principle. If the maximum principle does not hold, the revisions to the nonlinear potentials may lead to spurious solutions. Moreover, the second order scheme by stabilization only possesses the conditional energy stability (cf. the detailed rigorous proof in [38]), i.e., the time step size is controlled by the interfacial thickness.

For the BFS-PF model system, both of the two methods are not optimal choices. First, it is uncertain whether the solution of the system always satisfies certain maximum principle. Second, since the coupling entropy term involves two variables, it is particularly not clear and questionable whether the potential $\boldsymbol{E}_3(\phi, \rho)$ could be split into the combinations of a convex and a concave part. We aim to develop numerical schemes that could possess the advantages of both the convex splitting approach and the linear stabilization approach, but avoid their imperfections mentioned above. More precisely, we expect that the schemes are efficient (linear system), stable (unconditionally energy stable), and accurate (ready for second order or even higher order in time). To this end, we use a novel *Invariant Energy Quadratization* (IEQ) approach to design desired numerical schemes, without worrying about whether the continuous/discrete maximum principle holds or a convexity/concavity splitting exists.

Following the work in [7, 10], we first regularize the F-H potential from domain $(0,1)$ to $(-\infty,\infty)$, where the logarithmic functional is replaced by a $C^2$ continuous, convex, piecewise function. More precisely, for any $\widehat{\epsilon} > 0$, the regularized F-H potential is given by

$$\widehat{G}(\rho) = \begin{cases} \rho\ln\rho + \frac{(1-\rho)^2}{2\widehat{\epsilon}} + (1-\rho)\ln\widehat{\epsilon} - \frac{\widehat{\epsilon}}{2}, & \text{if} \quad \rho \geq 1-\widehat{\epsilon}, \\ \rho\ln\rho + (1-\rho)\ln(1-\rho), & \text{if} \quad \widehat{\epsilon} \leq \rho \leq 1-\widehat{\epsilon}, \\ (1-\rho)\ln(1-\rho) + \frac{\rho^2}{2\widehat{\epsilon}} + \rho\ln\widehat{\epsilon} - \frac{\widehat{\epsilon}}{2}, & \text{if} \quad \rho \leq \widehat{\epsilon}. \end{cases} \tag{3.1}$$

When $\widehat{\epsilon} \to 0$, $\widehat{G}(\rho) \to G(\rho)$. It was proved in [7, 10] that the error bound between the regularized system and the original system is controlled by $\widehat{\epsilon}$ up to a constant for the Cahn-Hilliard equation. Thus we consider the numerical solution to the model formulated with the regularized functional $\widehat{G}(\rho)$, but omit the $\widehat{}$ in the notation for convenience.

3.1. **Transformed system.** It is obvious that $G(\rho)$ is bounded from below although it is not always positive in the whole domain. Thus we could rewrite the free energy functional to the following equivalent form:

$$\boldsymbol{E}_{tot}(\phi,\rho) = \int_\Omega \Big(\frac{\epsilon}{2}|\nabla\phi|^2 + \frac{1}{4\epsilon}(\phi^2-1)^2 + \frac{\alpha}{2}(\rho - |\nabla\phi|)^2 + \beta\Big(\sqrt{G(\rho)+B}\,\Big)^2\Big)d\boldsymbol{x} - \beta B|\Omega|, \tag{3.2}$$

where $B$ is some positive constant to ensure $G(\rho) + B > 0$, for example, $B = 1$ that is adpoted in the simulations. We emphasize that the free energy is invariant because we simply add a zero term $\beta B - \beta B$ therein. Then we define three auxiliary variables to be the square root of the nonlinear potentials $F(\phi)$, $(\rho - |\nabla\phi|)^2$ and $G(\rho) + B$ by

$$U = \phi^2 - 1, \tag{3.3}$$

$$V = \rho - |\nabla\phi|, \tag{3.4}$$

$$W = \sqrt{G(\rho)+B}. \tag{3.5}$$

In turn, the total free energy (3.2) can be transformed as

$$\boldsymbol{E}_{tot}(\phi, U, V, W) = \int_\Omega \Big(\frac{\epsilon}{2}|\nabla\phi|^2 + \frac{1}{4\epsilon}U^2 + \frac{\alpha}{2}V^2 + \beta W^2\Big)d\boldsymbol{x} - \beta B|\Omega|. \tag{3.6}$$

Then we obtain a new, but equivalent partial differential system as follows:

$$\phi_t = M_1\Delta\mu_\phi, \tag{3.7}$$

$$\mu_\phi = -\epsilon\Delta\phi + \frac{1}{\epsilon}\phi U + \alpha\nabla\cdot(V\boldsymbol{Z}), \tag{3.8}$$

$$\rho_t = M_2\Delta\mu_\rho, \tag{3.9}$$

$$\mu_\rho = \alpha V + \beta H W, \tag{3.10}$$

$$U_t = 2\phi\phi_t, \tag{3.11}$$

$$V_t = \rho_t - \boldsymbol{Z}\cdot\nabla\phi_t, \tag{3.12}$$

$$W_t = \frac{1}{2}H\rho_t, \tag{3.13}$$

with $H = \frac{g(\rho)}{\sqrt{G(\rho)+B}}$ where $g(\rho) = G'(\rho)$.

The initial conditions are correspondingly

$$\begin{cases} \phi|_{(t=0)} = \phi_0, \quad \rho|_{(t=0)} = \rho_0, \\ U|_{(t=0)} = \phi_0^2 - 1,\ V|_{(t=0)} = \rho_0 - |\nabla\phi_0|,\ W|_{(t=0)} = \sqrt{G(\rho_0)+B}. \end{cases} \tag{3.14}$$

The transformed system still follows the energy dissipation law. By taking the sum of the $L^2$ inner products of (3.7) with $-\mu_\phi$, of (3.8) with $\phi_t$, of (3.9) with $-\mu_\rho$, of (3.10) with $\rho_t$, of (3.11) with $\frac{1}{2\epsilon}U$, of (3.12) with $\alpha V$, of (3.13) with $2\beta W$, we obtain the energy dissipation law of the new system as (3.7)-(3.13)

as

$$\frac{d}{dt}\boldsymbol{E}_{tot}(\phi, U, V, W) = -M_1\|\nabla\mu_\phi\|^2 - M_2\|\nabla\mu_\rho\|^2 \leq 0. \tag{3.15}$$

**Remark 3.1.** *We consider $F(\phi)$, $(\rho - |\nabla\phi|)^2$ and $G(\rho) + B$ as three quadratic functionals by applying appropriate substitutions if needed. Therefore, after simple substitutions using new variables $U, V, W$ defined in* (3.3)-(3.4)*, the energy is transformed to an equivalent quadratic form. We emphasize that the new transformed system* (3.7)–(3.13) *is exactly equivalent to the original system* (2.6)-(2.9) *since* (3.3)-(3.5) *can be easily obtained by integrating* (3.11)-(3.13) *with respect to the time. Therefore, the energy law* (3.15) *for the transformed system is exactly the same as the energy law* (2.10) *for the original system for the time-continuous case. We will develop time-marching schemes for the new transformed system* (3.7)–(3.13) *that in turn follow the new energy dissipation law* (3.15) *instead of the energy law for the original system* (2.10)*.*

In the following, we focus on designing numerical schemes for time stepping of the transformed system (3.7)-(3.13), that are linear and satisfy discrete analogues of the energy law (3.15). Let $\delta t > 0$ denote the time step size and set $t_n = n\,\delta t$ for $0 \leq n \leq N$ with the ending time $T = N\,\delta t$.

3.2. **First order scheme.** We now present a first order time stepping scheme to solve the system (3.7)-(3.13). Assuming that $\phi^n, \rho^n, U^n, V^n, W^n$ are already known, we then solve $\phi^{n+1}$, $\rho^{n+1}$, $U^{n+1}$, $V^{n+1}$ and $W^{n+1}$ from the following first order temporal semi-discretized system:

$$\frac{\phi^{n+1} - \phi^n}{\delta t} = M_1\Delta\mu_\phi^{n+1}, \tag{3.16}$$

$$\mu_\phi^{n+1} = -\epsilon\Delta\phi^{n+1} + \frac{1}{\epsilon}\phi^n U^{n+1} + \alpha\nabla\cdot(V^{n+1}\boldsymbol{Z}^n), \tag{3.17}$$

$$\frac{\rho^{n+1} - \rho^n}{\delta t} = M_2\Delta\mu_\rho^{n+1}, \tag{3.18}$$

$$\mu_\rho^{n+1} = \alpha V^{n+1} + \beta H^n W^{n+1}, \tag{3.19}$$

$$U^{n+1} - U^n = 2\phi^n(\phi^{n+1} - \phi^n), \tag{3.20}$$

$$V^{n+1} - V^n = (\rho^{n+1} - \rho^n) - \boldsymbol{Z}^n\cdot\nabla(\phi^{n+1} - \phi^n), \tag{3.21}$$

$$W^{n+1} - W^n = \frac{1}{2}H^n(\rho^{n+1} - \rho^n) \tag{3.22}$$

with periodic boundary condition being imposed and $\boldsymbol{Z}^n = \boldsymbol{Z}(\phi^n), H^n = H(\phi^n)$.

The impressing part of the above schemes is that all nonlinear coefficient of the new variables $U, V, W$ are treated explicitly, which can tremendously simply the computations. Moreover, in fact, we can rewrite the equations (3.20)-(3.22) as follows,

$$\begin{cases} U^{n+1} = A_1 + 2\phi^n\phi^{n+1}, \\ V^{n+1} = B_1 + \rho^{n+1} - \boldsymbol{Z}^n\cdot\nabla\phi^{n+1}, \\ W^{n+1} = C_1 + \dfrac{1}{2}H^n\rho^{n+1}, \end{cases} \tag{3.23}$$

where

$$\begin{cases} A_1 = U^n - 2(\phi^n)^2, \\ B_1 = V^n - \rho^n + \boldsymbol{Z}^n\cdot\nabla\phi^n, \\ C_1 = W^n - \dfrac{1}{2}H^n\rho^n. \end{cases} \tag{3.24}$$

Thus the system (3.16)-(3.22) can be rewritten as

$$\frac{\phi^{n+1}-\phi^n}{\delta t} = M_1\Delta\mu_\phi^{n+1}, \tag{3.25}$$

$$\mu_\phi^{n+1} = -\epsilon\Delta\phi^{n+1} + P_1(\phi^{n+1},\rho^{n+1}) + S_1^n, \tag{3.26}$$

$$\frac{\rho^{n+1}-\rho^n}{\delta t} = M_2\Delta\mu_\rho^{n+1}, \tag{3.27}$$

$$\mu_\rho^{n+1} = Q_1(\phi^{n+1},\rho^{n+1}) + S_2^n, \tag{3.28}$$

where

$$\begin{cases} P_1(\phi,\rho) = \frac{1}{\epsilon}2\phi^n\phi^n\phi + \alpha\nabla\cdot(\rho\boldsymbol{Z}^n) - \alpha\nabla\cdot((\boldsymbol{Z}^n\cdot\nabla\phi)\boldsymbol{Z}^n),\\ Q_1(\phi,\rho) = \alpha\rho - \alpha\boldsymbol{Z}^n\cdot\nabla\phi + \frac{1}{2}\beta H^nH^n\rho,\\ S_1^n = \frac{1}{\epsilon}\phi^nA_1 + \alpha\nabla\cdot(B_1\boldsymbol{Z}^n),\\ S_2^n = \alpha B_1 + \beta H^nC_1. \end{cases} \tag{3.29}$$

Therefore, we can solve $(\phi,\mu_\phi,\rho,\mu_\rho)^{n+1}$ directly from (3.25)–(3.28). Once we obtain $\phi^{n+1},\rho^{n+1}$, the $U^{n+1},V^{n+1},W^{n+1}$ are automatically given in (3.20)-(3.22). Namely, the new variables $U,V,W$ do not involve any extra computational costs.

Meanwhile, we notice the following two facts hold for linear operators $P_1(\phi,\rho)$ and $Q_1(\phi,\rho)$.

- If $X_1,X_2,Y_1,Y_2$ satisfy periodic boundary conditions, we have

$$(P_1(X_1,X_2),Y_1) + (Q_1(X_1,X_2),Y_2) = (P_1(Y_1,Y_2),X_1) + (Q_1(Y_1,Y_2),X_2). \tag{3.30}$$

- For any $X_1,X_2$ with $\int_\Omega X_1d\boldsymbol{x} = \int_\Omega X_2d\boldsymbol{x} = 0$, we have

$$(P_1(X_1,X_2),X_1) + (Q_1(X_1,X_2),X_2) = \frac{2}{\epsilon}\|\phi^nX_1\|^2 + \frac{1}{2}\beta\|H^nX_2\|^2 + \alpha\|X_2 - \boldsymbol{Z}^n\cdot\nabla X_1\|^2 \geq 0, \tag{3.31}$$

  where "=" is valid if and only if $X_1 = X_2 = 0$ pointwise assuming that $\phi^n,H^n,Z^n$ are not zero pointwise (if $\phi^n,H^n,Z^n$ are all zeros, then the equations only have zero solutions).

The first theorem ensures the efficiency of the linear scheme (3.25)-(3.28) as follows.

**Theorem 3.1.** *The linear system of* (3.25)-(3.28) *is symmetric (self-adjoint) and positive definite for the variable* $\rho^{n+1},\phi^{n+1}$.

*Proof.* From (3.25) and (3.27), by taking the $L^2$ inner product with 1, we derive

$$\begin{aligned} \int_\Omega\rho^{n+1}d\boldsymbol{x} &= \int_\Omega\rho^nd\boldsymbol{x} = \cdots = \int_\Omega\rho^0d\boldsymbol{x},\\ \int_\Omega\phi^{n+1}d\boldsymbol{x} &= \int_\Omega\phi^nd\boldsymbol{x} = \cdots = \int_\Omega\phi^0d\boldsymbol{x}. \end{aligned} \tag{3.32}$$

Let

$$\begin{aligned} \alpha_\rho^0 &= \frac{1}{|\Omega|}\int_\Omega\rho^0d\boldsymbol{x}, & \alpha_\phi^0 &= \frac{1}{|\Omega|}\int_\Omega\phi^0d\boldsymbol{x},\\ \beta_\rho^\mu &= \frac{1}{|\Omega|}\int_\Omega\mu_\rho^{n+1}d\boldsymbol{x}, & \beta_\phi^\mu &= \frac{1}{|\Omega|}\int_\Omega\mu_\phi^{n+1}d\boldsymbol{x}, \end{aligned} \tag{3.33}$$

and we define

$$\begin{aligned} \widehat{\rho}^{n+1} &= \rho^{n+1} - \alpha_\rho^0, & \widehat{\mu}_\rho^{n+1} &= \mu_\rho^{n+1} - \beta_\rho^\mu.\\ \widehat{\phi}^{n+1} &= \phi^{n+1} - \alpha_\phi^0, & \widehat{\mu}_\phi^{n+1} &= \mu_\phi^{n+1} - \beta_\phi^\mu. \end{aligned} \tag{3.34}$$

Thus from (3.25)-(3.28), $\widehat{\rho}^{n+1}, \widehat{\phi}^{n+1}, \widehat{\mu}_\rho^{n+1}, \widehat{\mu}_\phi^{n+1}$ are the solutions for the following equations

$$\frac{\phi}{\delta t} - M_1\Delta\mu_\phi = \frac{\widehat{\phi}^n}{\delta t}, \tag{3.35}$$

$$\mu_\phi + \beta_\phi^\mu + \epsilon\Delta\phi - P_1(\phi,\rho) = f_1^n \tag{3.36}$$

$$\frac{\rho}{\delta t} - M_1\Delta\mu_\rho = \frac{\widehat{\rho}^n}{\delta t}, \tag{3.37}$$

$$\mu_\rho + \beta_\rho^\mu - Q_1(\phi,\rho) = f_2^n. \tag{3.38}$$

with $\int_\Omega \phi d\boldsymbol{x} = \int_\Omega \rho d\boldsymbol{x} = 0$.

We define the inverse Laplace operator $u$ (with $\int_\Omega u d\boldsymbol{x} = 0$) $\to v := \Delta^{-1}u$ by

$$\begin{cases} \Delta v = u, \\ \displaystyle\int_\Omega v d\boldsymbol{x} = 0, \end{cases} \tag{3.39}$$

with the periodic boundary conditions.

Applying $-\Delta^{-1}$ to (3.35), (3.37) and using (3.36), (3.38), we obtain

$$-\frac{1}{M_1\delta t}\Delta^{-1}\phi - \epsilon\Delta\phi + P_1(\phi,\rho) - \beta_\phi^\mu = -\frac{1}{M_1\delta t}\Delta^{-1}\frac{\widehat{\phi}^n}{\delta t} - f_1^n, \tag{3.40}$$

$$-\frac{1}{M_2\delta t}\Delta^{-1}\rho + Q_1(\phi,\rho) - \beta_\rho^\mu = -\frac{1}{M_2\delta t}\Delta^{-1}\frac{\widehat{\rho}^n}{\delta t} - f_2^n. \tag{3.41}$$

We denote the above linear system as $\mathbb{A}\boldsymbol{X} = \mathbf{b}$ where $\boldsymbol{X} = [X_1, X_2]^T$ where $\phi, \rho$ are denoted by $X_1, X_2$ respectively.

For any $\boldsymbol{X} = [X_1, X_2]^T$ and $\boldsymbol{Y} = [Y_1, Y_2]^T$ with $\int_\Omega X_1 d\boldsymbol{x} = \int_\Omega X_2 d\boldsymbol{x} = \int_\Omega Y_1 d\boldsymbol{x} = \int_\Omega Y_2 d\boldsymbol{x} = 0$, from (3.30), we derive

$$(\mathbb{A}\boldsymbol{X}, \boldsymbol{Y}) = (\mathbb{A}\boldsymbol{Y}, \boldsymbol{X}). \tag{3.42}$$

Therefore $\mathbb{A}$ is self adjoint. Meanwhile, we derive

$$\begin{aligned}(\mathbb{A}\boldsymbol{X}, \boldsymbol{X}) =& (-\frac{1}{M_1\delta t}\Delta^{-1}X_1, X_1) + \epsilon(\nabla X_1, \nabla X_1) + (P_1(X_1,X_2), X_1) - (\beta_\phi^\mu, X_1) \\ &+ (-\frac{1}{M_2\delta t}\Delta^{-1}X_2, X_2) + (Q_1(X_1,X_2), X_2) - (\beta_\rho^\mu, X_2) \\ =& \frac{1}{M_1\delta t}\|X_1\|_{H^{-1}}^2 + \frac{1}{M_2\delta t}\|X_2\|_{H^{-1}}^2 \\ &+ \epsilon\|\nabla X_1\|^2 + \frac{2}{\epsilon}\|\phi^n X_1\|^2 + \frac{\beta}{2}\|H^n X_2\|^2 + \alpha\|X_2 - \boldsymbol{Z}^n \cdot \nabla X_1\|^2 \geq 0\end{aligned} \tag{3.43}$$

where "=" is valid if and only if $X_1 = X_2 = 0$ from the fact of (3.31). Thus we conclude the theorem. □

**Remark 3.2.** *We can show the well-posedness of the linear system* $\mathbb{A}\boldsymbol{X} = \mathbf{b}$ *from the Lax-Milgram theorem. Define* $(\boldsymbol{X}, \boldsymbol{Y})_{\mathbb{A}} = (\mathbb{A}\boldsymbol{X}, \boldsymbol{Y})$, $\|\boldsymbol{X}\|_{\mathbb{A}} = \sqrt{(\boldsymbol{X}, \boldsymbol{X})_{\mathbb{A}}}$ *for any* $\boldsymbol{X}, \boldsymbol{Y} \in L^2_{per}(\Omega)$ *and the subset* $\mathbb{S} = \{\boldsymbol{X} \in L^2_{per}(\Omega) : \|\boldsymbol{X}\|_{\mathbb{A}} < \infty\}$. *It is not hard to show the* $\|\cdot\|_{\mathbb{A}}$ *is a well defined norm and* $\mathbb{S}$ *is a Hilbert subspace of* $L^2_{per}(\Omega)$ *associated with the inner product* $(\cdot,\cdot)_{\mathbb{A}}$. *The details of the derivation process are omitted since the proof only involves the standard analysis techniques. Therefore the well-posedness of* $\mathbb{A}\boldsymbol{X} = \mathbf{b}$ *in the weak sense is valid from the Lax-Milgram theorem, i.e., the system admits a unique weak solution in* $\mathbb{S}$.

The stability result of the proposed first order scheme is given below that follows the same lines as in the derivation of the new energy dissipation law (3.15).

**Theorem 3.2.** *The first order linear scheme* (3.16)-(3.22) *is unconditionally energy stable, i.e., satisfies the following discrete energy dissipation law:*

$$\frac{1}{\delta t}(E_{1st}^{n+1} - E_{1st}^n) \leq -M_1\|\nabla\mu_\phi^{n+1}\|^2 - M_2\|\nabla\mu_\rho^{n+1}\|^2, \tag{3.44}$$

*where*

$$E^n_{1st} = \frac{\epsilon}{2}\|\nabla\phi^n\|^2 + \frac{1}{4\epsilon}\|U^n\|^2 + \frac{\alpha}{2}\|V^n\|^2 + \beta\|W^n\|^2 - \beta B. \tag{3.45}$$

*Proof.* By taking the $L^2$ inner product of (3.16) with $-\mu_\phi^{n+1}$, we obtain

$$\Big(\frac{\phi^{n+1}-\phi^n}{\delta t}, -\mu_\phi^{n+1}\Big) = M_1\|\nabla\mu_\phi^{n+1}\|^2. \tag{3.46}$$

By taking the $L^2$ inner product of (3.17) with $\frac{\phi^{n+1}-\phi^n}{\delta t}$ and applying the following identities

$$2(a-b,a) = |a|^2 - |b|^2 + |a-b|^2, \tag{3.47}$$

we obtain

$$\begin{aligned}\Big(\mu_\phi^{n+1}, \frac{\phi^{n+1}-\phi^n}{\delta t}\Big) =& \frac{1}{\delta t}\frac{\epsilon}{2}(\|\nabla\phi^{n+1}\|^2 - \|\nabla\phi^{n+1}\|^2 + \|\nabla(\phi^{n+1}-\phi^n)\|^2)\\ &+ \frac{1}{\epsilon}\Big(\phi^n U^{n+1}, \frac{\phi^{n+1}-\phi^n}{\delta t}\Big) - \alpha\Big(V^{n+1}Z(\phi^n), \nabla\frac{\phi^{n+1}-\phi^n}{\delta t}\Big).\end{aligned} \tag{3.48}$$

By taking the $L^2$ inner product of (3.18) with $-\mu_\rho^{n+1}$, we obtain

$$-\Big(\frac{\rho^{n+1}-\rho^n}{\delta t}, \mu_\rho^{n+1}\Big) = M_2\|\nabla\mu_\rho^{n+1}\|^2. \tag{3.49}$$

By taking the $L^2$ inner product of (3.19) with $\frac{\rho^{n+1}-\rho^n}{\delta t}$, we obtain

$$\Big(\mu_\rho^{n+1}, \frac{\rho^{n+1}-\rho^n}{\delta t}\Big) = \alpha\Big(V^{n+1}, \frac{\rho^{n+1}-\rho^n}{\delta t}\Big) + \beta\Big(H(\rho^n)W^{n+1}, \frac{\rho^{n+1}-\rho^n}{\delta t}\Big). \tag{3.50}$$

By taking the $L^2$ inner product of (3.20) with $\frac{1}{2\epsilon\delta t}U^{n+1}$, we obtain

$$\frac{1}{4\delta t\epsilon}\Big(\|U^{n+1}\|^2 - \|U^n\|^2 + \|U^{n+1}-U^n\|^2\Big) = \frac{1}{\epsilon}\Big(\phi^n\frac{\phi^{n+1}-\phi^n}{\delta t}, U^{n+1}\Big). \tag{3.51}$$

By taking the $L^2$ inner product of (3.21) with $\frac{1}{\delta t}\alpha V^{n+1}$, we obtain

$$\begin{aligned}&\frac{\alpha}{2\delta t}\Big(\|V^{n+1}\|^2 - \|V^n\|^2 + \|V^{n+1}-V^n\|^2\Big)\\ &= \alpha\Big(\frac{\rho^{n+1}-\rho^n}{\delta t}, V^{n+1}\Big) - \alpha\Big(Z(\phi^n)\nabla\frac{\phi^{n+1}-\phi^n}{\delta t}, V^{n+1}\Big).\end{aligned} \tag{3.52}$$

By taking the $L^2$ inner product of (3.22) with $2\frac{1}{\delta t}\beta W^{n+1}$, we obtain

$$\frac{\beta}{\delta t}\Big(\|W^{n+1}\|^2 - \|W^n\|^2 + \|W^{n+1}-W^n\|^2\Big) =\beta\Big(H(\rho^n)\frac{\rho^{n+1}-\rho^n}{\delta t}, W^{n+1}\Big). \tag{3.53}$$

Combination of (3.46), (3.48)-(3.53) gives us

$$\begin{aligned}&\frac{\epsilon}{2}(\|\nabla\phi^{n+1}\|^2 - \|\nabla\phi^{n+1}\|^2 + \|\phi^{n+1}-\phi^n\|^2) + \frac{1}{4\epsilon}(\|U^{n+1}\|^2 - \|U^n\|^2 + \|U^{n+1}-U^n\|^2)\\ &+ \frac{\alpha}{2}(\|V^{n+1}\|^2 - \|V^n\|^2 + \|U^{n+1}-V^n\|^2) + \beta(\|W^{n+1}\|^2 - \|W^n\|^2 + \|W^{n+1}-W^n\|^2)\\ &= -M_1\delta t\|\nabla\mu_\phi^{n+1}\|^2 - M_2\delta t\|\nabla\mu_\rho^{n+1}\|^2.\end{aligned} \tag{3.54}$$

Finally, we obtain the desired result (3.44) after dropping some positive terms from (3.54). □

The proposed scheme follows the new energy dissipation law (3.15) instead of the energy law for the originated system (2.10). For time-continuous case, (3.15) and (2.10) are identical. For time-discrete case, the discrete energy law (3.45) is the first order approximation to the new energy law (3.15). Moreover, the discrete energy functional $E(\phi^{n+1}, U^{n+1}, V^{n+1}, W^{n+1})$ is also the first order approximation to $E(\phi^{n+1}, \rho^{n+1})$ (defined in (2.5)), since $U^{n+1}, V^{n+1}, W^{n+1}$ are the first order approximations to $\phi^2 - 1$, $\rho - |\nabla\phi|$ and $\sqrt{G(\rho) + B)}$ respectively, that can be observed from the following facts, heuristically. We rewrite (3.20) as follows,

$$U^{n+1}-((\phi^{n+1})^2-1)=U^n-((\phi^n)^2-1)+R_{n+1}, \tag{3.55}$$

where

$$R_{n+1}=(\phi^{n+1}-\phi^n)^2. \tag{3.56}$$

Since $R_{n+1}$ is the resudue of the Taylor expansion for the function $(\phi^{n+1})^2-1$ at $t=t^n$, thus $R_{n+1}=O(\delta t^2)$ for $0\le k\le n+1$. Note that $U^0=(\phi^0)^2-1$, by mathematical induction we can easily get $U^{n+1}=(\phi^{n+1})^2-1+O(\delta t)$. The same argument can be performed for the other variable $V^{n+1}$ and $W^{n+1}$ since (3.21) and (3.22) are the Tylor expansion formulation for $V^{n+1}$ and $W^{n+1}$ at $t=t^n$, respectively.

The energy stability of the numerical schemes is formally derived and it is expected that the optimal error estimates can be obtained correspondingly. This expectation is supported by the error analysis of the IEQ scheme for the classical Cahn-Hilliard equation with the double well potential, that is quite clear if we notice that the $H^1$ bound for the numerical solution is obtained due to the Poincare inequality. For the BFS-PF model with multiple variables, the error analysis is more complicated and we will implement it in the future work that will follow the same lines as the literatures [12, 38, 48].

**Remark 3.3.** *We note that the idea of the IEQ approach is very simple but quite different from the traditional time marching schemes. For example, it does not require the convexity as the convex splitting approach (cf. [11]) or the boundedness for the second order derivative as the stabilization approach (cf. [6, 38, 49]). Through a simple substitution of new variables, the complicated nonlinear potentials are transformed into quadratic forms. We summarize the great advantages of this quadratic transformations as follows: (i) this quadratization method works well for various complex nonlinear terms as long as the corresponding nonlinear potentials are bounded from below; (ii) the complicated nonlinear potential is transferred to a quadratic polynomial form which is much easier to handle; (iii) the derivative of the quadratic polynomial is linear, which provides the fundamental support for linearization method; (iv) the quadratic formulation in terms of new variables can automatically maintain this property of positivity (or bounded from below) of the nonlinear potentials.*

*When the nonlinear potential is the fourth order polynomial, e.g., the double well potential, the IEQ method is exactly the same as the so-called Lagrange Multiplier method developed in [16]. We remark that the idea of Lagrange Multiplier method only works well for the fourth order polynomial potential ($\phi^4$). This is because the nonlinear term $\phi^3$ (the derivative of $\phi^4$) can be naturally decomposed into a multiplication of two factors: $\lambda(\phi)\phi$ that is the Lagrange multiplier term, and the $\lambda(\phi)=\phi^2$ is then defined as the new auxiliary variable $U$. However, this method might not succeed for other type potentials. For instance, for the F-H type potential where the nonlinear term is $\ln(\frac{\phi}{1-\phi})$, if one forcefully rewrites this term as $\lambda(\phi)\phi$, then $\lambda(\phi)=\frac{\ln(\frac{\phi}{1-\phi})}{\phi}$ that is the format for the new variable $U$. Obviously, such a form is unworkable for algorithms design. About the application of the IEQ approach to handle other type of nonlinear potentials, we refer to the authors' other work in [51, 59].*

**Remark 3.4.** *The IEQ approach provides more efficiency than the nonlinear approach. More precisely, if the nonlinear potential takes the form of a convex polynomial, i.e., $E(\phi)=\int_\Omega \phi^{2K}d\boldsymbol{x}$ ($K\ge 2$), then we will have the linear scheme that includes $(\phi^n)^{2K-2}\phi^{n+1}$. Let us consider the implicit or convex splitting approach where the nonlinear term has to be discretized by $(\phi^{n+1})^{2K-1}$. Therefore, if the Newton iterative method is applied for this term, at each iteration the nonlinear convex splitting approach would yield the same linear operator as IEQ approach. Hence the cost of solving the IEQ scheme is the same as the cost of performing one iteration of Newton method for the implicit/convex splitting approach, provided that the same linear solvers are applied (for instance multi-grid with Gauss-Seidel relaxation). It is clear that the IEQ scheme would be much more efficient than the nonlinear schemes.*

3.3. **Second order scheme.** Undoubtedly, higher order time marching schemes are preferable to lower order schemes if the adopted time step is expected to be as large as possible under certain accuracy requests. Thus we construct the second order time stepping schemes to solve the system (3.7)-(3.13), based on the Adam-Bashforth backward differentiation formulas (BDF2).

Assume that $\phi^{n-1}$, $\rho^{n-1}$, $U^{n-1}$, $V^{n-1}$, $W^{n-1}$ and $\phi^n$, $\rho^n$, $U^n$, $V^n$, $W^n$ are known, we solve $\phi^{n+1}$, $\rho^{n+1}$, $U^{n+1}$, $V^{n+1}$, $W^{n+1}$ as follows:

$$\frac{3\phi^{n+1}-4\phi^n+\phi^{n-1}}{2\delta t} = M_1\Delta\mu_\phi^{n+1}, \tag{3.57}$$

$$\mu_\phi^{n+1} = -\epsilon\Delta\phi^{n+1} + \frac{1}{\epsilon}\phi^\star U^{n+1} + \alpha\nabla\cdot(V^{n+1}\boldsymbol{Z}^\star), \tag{3.58}$$

$$\frac{3\rho^{n+1}-4\rho^n+\rho^{n-1}}{2\delta t} = M_2\Delta\mu_\rho^{n+1}, \tag{3.59}$$

$$\mu_\rho^{n+1} = \alpha V^{n+1} + \beta H^\star W^{n+1}, \tag{3.60}$$

$$3U^{n+1} - 4U^n + U^{n-1} = 2\phi^\star(3\phi^{n+1}-4\phi^n+\phi^{n-1}), \tag{3.61}$$

$$3V^{n+1} - 4V^n + V^{n-1} = (3\rho^{n+1}-4\rho^n+\rho^{n-1}) - \boldsymbol{Z}^\star\cdot\nabla(3\phi^{n+1}-4\phi^n+\phi^{n-1}), \tag{3.62}$$

$$3W^{n+1} - 4W^n + W^{n-1} = \frac{1}{2}H^\star(3\rho^{n+1}-4\rho^n+\rho^{n-1}), \tag{3.63}$$

with periodic boundary condition being imposed, where

$$\phi^\star = 2\phi^n - \phi^{n-1}, \boldsymbol{Z}^\star = \boldsymbol{Z}(2\phi^n-\phi^{n-1}), H^\star = H(2\rho^n-\rho^{n-1}). \tag{3.64}$$

Similar to the first order scheme, we can rewrite the equations (3.61)–(3.63) as follows,

$$\begin{cases} U^{n+1} = A_2 + 2\phi^\star\phi^{n+1}, \\ V^{n+1} = B_2 + \rho^{n+1} - \boldsymbol{Z}^\star\cdot\nabla\phi^{n+1}, \\ W^{n+1} = C_2 + \frac{1}{2}H^\star\rho^{n+1}, \end{cases} \tag{3.65}$$

where

$$\begin{cases} A_2 = U^\dagger - 2\phi^\star\phi^\dagger, \\ B_2 = V^\dagger - \rho^\dagger + \boldsymbol{Z}^\star\cdot\nabla\phi^\dagger, \\ C_2 = W^\dagger - \frac{1}{2}H^\star\rho^\dagger, \end{cases} \tag{3.66}$$

where $S^\dagger = \frac{4S^n - S^{n-1}}{3}$ for any variable $S$. Thus the system (3.57)-(3.63) can be rewritten as

$$\frac{3\phi^{n+1}-4\phi^n+\phi^{n-1}}{2\delta t} = M_1\Delta\mu_\phi^{n+1}, \tag{3.67}$$

$$\mu_\phi^{n+1} = -\epsilon\Delta\phi^{n+1} + P_2(\phi^{n+1},\rho^{n+1}) + R_1^n, \tag{3.68}$$

$$\frac{3\rho^{n+1}-4\rho^n+\rho^{n-1}}{2\delta t} = M_2\Delta\mu_\rho^{n+1}, \tag{3.69}$$

$$\mu_\rho^{n+1} = Q_2(\phi^{n+1},\rho^{n+1}) + R_2^n, \tag{3.70}$$

where

$$\begin{cases} P_2(\phi,\rho) = \frac{1}{\epsilon}2\phi^\star\phi^\star\phi + \alpha\nabla\cdot(\rho\boldsymbol{Z}^\star) - \alpha\nabla\cdot((\boldsymbol{Z}^\star\cdot\nabla\phi)\boldsymbol{Z}^\star), \\ Q_2(\phi,\rho) = \alpha\rho - \alpha\boldsymbol{Z}^\star\cdot\nabla\phi + \frac{1}{2}\beta H^\star H^\star\rho, \\ R_1^n = \frac{1}{\epsilon}\phi^\star A_2 + \alpha\nabla\cdot(B_2\boldsymbol{Z}^\star), \\ R_2^n = \alpha B_2 + \beta H^\star C_2. \end{cases} \tag{3.71}$$

For $P_2$ and $Q_2$, we derive similar properties as $P_1, Q_1$.

- If $X_1, X_2, Y_1, Y_2$ satisfy periodic boundary conditions, we have

$$(P_2(X_1,X_2),Y_1) + (Q_2(X_1,X_2),Y_2) = (P_2(Y_1,Y_2),X_1) + (Q_2(Y_1,Y_2),X_2), \tag{3.72}$$

- For any $X_1, X_2$ with $\int_\Omega X_1 d\boldsymbol{x} = \int_\Omega X_2 d\boldsymbol{x} = 0$, we have

$$(P_2(X_1,X_2),X_1)+(Q_2(X_1,X_2),X_2)=\frac{2}{\epsilon}\|\phi^\star X_1\|^2+\frac{1}{2}\beta\|H^\star X_2\|^2+\alpha\|X_2-\boldsymbol{Z}^\star\cdot\nabla X_1\|^2\geq 0, \tag{3.73}$$

where "=" is valid if and only if $X_1 = X_2 = 0$ pointwise.

**Theorem 3.3.** *The linear system of* (3.67)-(3.70) *is symmetric positive definite for the variable* $\phi^{n+1}, \rho^{n+1}$.

*Proof.* We omit the details here since the proof is similar to that for Theorem 3.1. □

**Theorem 3.4.** *The second order linear scheme* (3.57)–(3.63) *is unconditionally energy stable, i.e., satisfies the following discrete energy dissipation law:*

$$\frac{1}{\delta t}(E_{2nd}^{n+1,n}-E_{2nd}^{n,n-1})\leq -M_1\|\nabla\mu_\phi^{n+1}\|^2-M_2\|\nabla\mu_\rho^{n+1}\|^2, \tag{3.74}$$

*where*

$$\begin{aligned}E_{2nd}^{n+1,n}&=\frac{\epsilon}{2}\Big(\frac{1}{2}\|\nabla\phi^{n+1}\|^2+\frac{1}{2}\|2\nabla\phi^{n+1}-\nabla\phi^n\|^2\Big)+\frac{1}{4\epsilon}\Big(\frac{1}{2}\|U^{n+1}\|^2+\frac{1}{2}\|2U^{n+1}-U^n\|^2\Big)\\&+\frac{\alpha}{2}\Big(\frac{1}{2}\|V^{n+1}\|^2+\frac{1}{2}\|2V^{n+1}-V^n\|^2\Big)+\beta\Big(\frac{1}{2}\|W^{n+1}\|^2+\frac{1}{2}\|2W^{n+1}-W^n\|^2\Big).\end{aligned}$$

*Proof.* By taking the $L^2$ inner product of (3.57) with $-\mu_\phi^{n+1}$, we obtain

$$-\Big(\frac{3\phi^{n+1}-4\phi^n+\phi^{n-1}}{2\delta t},\mu_\phi^{n+1}\Big)=M_1\|\nabla\mu_\phi^{n+1}\|^2. \tag{3.75}$$

By taking the $L^2$ inner product of (3.58) with $\frac{3\phi^{n+1}-4\phi^n+\phi^{n-1}}{2\delta t}$ and applying the following identities the following identity

$$2(3a-4b+c,a)=|a|^2-|b|^2+|2a-b|^2-|2b-c|^2+|a-2b+c|^2, \tag{3.76}$$

we obtain

$$\begin{aligned}&\Big(\mu_\phi^{n+1},\frac{3\phi^{n+1}-4\phi^n+\phi^{n-1}}{2\delta t}\Big)\\&=\frac{1}{2\delta t}\frac{\epsilon}{2}\Big(\|\nabla\phi^{n+1}\|^2-\|\nabla\phi^{n+1}\|^2+\|2\nabla\phi^{n+1}-\nabla\phi^n\|^2-\|2\nabla\phi^n-\nabla\phi^{n-1}\|^2\\&+\|\nabla\phi^{n+1}-2\nabla\phi^n+\nabla\phi^{n-1}\|^2\Big)+\frac{1}{\epsilon}\Big(\phi^\star U^{n+1},\frac{3\phi^{n+1}-4\phi^n+\phi^{n-1}}{2\delta t}\Big)\\&-\alpha\Big(V^{n+1}\boldsymbol{Z}^\star,\nabla\frac{3\phi^{n+1}-4\phi^n+\phi^{n-1}}{2\delta t}\Big).\end{aligned} \tag{3.77}$$

By taking the $L^2$ inner product of (3.59) with $-\mu_\rho^{n+1}$, we obtain

$$-\Big(\frac{3\rho^{n+1}-4\rho^n+\rho^{n-1}}{2\delta t},\mu_\rho^{n+1}\Big)=M_2\|\nabla\mu_\rho^{n+1}\|^2. \tag{3.78}$$

By taking the $L^2$ inner product of (3.60) with $\frac{3\rho^{n+1}-4\rho^n+\rho^{n-1}}{2\delta t}$, we obtain

$$\begin{aligned}\Big(\mu_\rho^{n+1},\frac{3\rho^{n+1}-4\rho^n+\rho^{n-1}}{2\delta t}\Big)&=\alpha\Big(V^{n+1},\frac{3\rho^{n+1}-4\rho^n+\rho^{n-1}}{2\delta t}\Big)\\&+\beta\Big(H^\star W^{n+1},\frac{3\rho^{n+1}-4\rho^n+\rho^{n-1}}{2\delta t}\Big).\end{aligned} \tag{3.79}$$

By taking the $L^2$ inner product of (3.61) with $\frac{1}{4\delta t\epsilon}U^{n+1}$, we obtain

$$\begin{aligned}&\frac{1}{8\delta t\epsilon}\Big(\|U^{n+1}\|^2-\|U^n\|^2+\|2U^{n+1}-U^n\|^2-\|2U^n-U^{n-1}\|^2+\|U^{n+1}-2U^n+U^{n-1}\|^2\Big)\\&=\frac{1}{\epsilon}\Big(\phi^\star\frac{3\phi^{n+1}-4\phi^n+\phi^{n-1}}{2\delta t},U^{n+1}\Big).\end{aligned} \tag{3.80}$$

By taking the $L^2$ inner product of (3.62) with $\frac{1}{2\delta t}\alpha V^{n+1}$, we obtain

$$\begin{aligned}&\frac{\alpha}{4\delta t}\Big(\|V^{n+1}\|^2-\|V^n\|^2+\|2V^{n+1}-V^n\|^2-\|2V^{n+1}-V^n\|^2+\|V^{n+1}-2V^n+V^{n-1}\|^2\Big)\\&\quad=\alpha\Big(\frac{3\rho^{n+1}-4\rho^n+\rho^{n-1}}{2\delta t},V^{n+1}\Big)-\alpha\Big(\boldsymbol{Z}^\star\nabla\frac{3\phi^{n+1}-4\phi^n+\phi^{n-1}}{2\delta t},V^{n+1}\Big).\end{aligned}\tag{3.81}$$

By taking the $L^2$ inner product of (3.63) with $\frac{\beta}{\delta t}W^{n+1}$, we obtain

$$\begin{aligned}&\frac{\beta}{2\delta t}\Big(\|W^{n+1}\|^2-\|W^n\|^2+\|2W^{n+1}-W^n\|^2-\|2W^{n+1}-W^n\|^2\\&\qquad+\|W^{n+1}-2W^n+W^{n-1}\|^2\Big)=\beta\Big(H^\star\frac{3\rho^{n+1}-4\rho^n+\rho^{n-1}}{2\delta t},W^{n+1}\Big).\end{aligned}\tag{3.82}$$

Combination of (3.75) and (3.77)–(3.82) gives

$$\begin{aligned}&\frac{\epsilon}{4}\Big(\|\nabla\phi^{n+1}\|^2-\|\nabla\phi^{n+1}\|^2+\|2\nabla\phi^{n+1}-\nabla\phi^n\|^2-\|2\nabla\phi^n-\nabla\phi^{n-1}\|^2\\&+\|\nabla\phi^{n+1}-2\nabla\phi^n+\nabla\phi^{n-1}\|^2\Big)\\&+\frac{1}{8\epsilon}\Big(\|U^{n+1}\|^2-\|U^n\|^2+\|2U^{n+1}-U^n\|^2-\|2U^n-U^{n-1}\|^2+\|U^{n+1}-2U^n+U^{n-1}\|^2\Big)\\&+\frac{\alpha}{4}\Big(\|V^{n+1}\|^2-\|V^n\|^2+\|2V^{n+1}-V^n\|^2-\|2V^{n+1}-V^n\|^2+\|V^{n+1}-2V^n+V^{n-1}\|^2\Big)\\&+\frac{\beta}{2}\Big(\|W^{n+1}\|^2-\|W^n\|^2+\|2W^{n+1}-W^n\|^2-\|2W^{n+1}-W^n\|^2+\|W^{n+1}-2W^n+W^{n-1}\|^2\Big)\\&=-M_1\delta t\|\nabla\mu_\phi^{n+1}\|^2-M_2\delta t\|\nabla\mu_\rho^{n+1}\|^2.\end{aligned}$$

Finally, we obtain the result (3.74) after dropping some positive terms from the above equation. □

Heuristically, the discrete energy law (3.74) is a second order approximation of $\frac{d}{dt}\boldsymbol{E}_{tot}$ in (3.15) since for any variable $S$, we have

$$\begin{aligned}&\Big(\frac{\|S^{n+1}\|^2+\|2S^{n+1}-S^n\|^2}{2\delta t}\Big)-\Big(\frac{\|S^n\|^2+\|2S^n-S^{n-1}\|^2}{2\delta t}\Big)\\&\qquad\cong\Big(\frac{\|S^{n+2}\|^2-\|S^n\|^2}{2\delta t}\Big)+o(\delta t^2)\cong\frac{d}{dt}\|S(t^{n+1})\|^2+o(\delta t^2).\end{aligned}$$

**Remark 3.5.** *One can easily develop the second order version based on the Crank-Nicolson type scheme, as follows, Assume that $\phi^{n-1}$, $\rho^{n-1}$, $U^{n-1}$, $V^{n-1}$, $W^{n-1}$ and $\phi^n$, $\rho^n$, $U^n$, $V^n$, $W^n$ are known, we solve $\phi^{n+1}$, $\rho^{n+1}$, $U^{n+1}$, $V^{n+1}$, $W^{n+1}$ as follows:*

$$\frac{\phi^{n+1}-\phi^n}{\delta t}=M_1\Delta\mu_\phi^{n+\frac12},\tag{3.83}$$

$$\mu_\phi^{n+\frac12}=-\epsilon\Delta\frac{\phi^{n+1}+\phi^n}{2}+\frac{1}{\epsilon}\phi^\circ\frac{U^{n+1}+U^n}{2}+\alpha\nabla\cdot(\frac{V^{n+1}+V^n}{2}\boldsymbol{Z}^\circ),\tag{3.84}$$

$$\frac{\rho^{n+1}-\rho^n}{\delta t}=M_2\Delta\mu_\rho^{n+\frac12},\tag{3.85}$$

$$\mu_\rho^{n+\frac12}=\alpha\frac{V^{n+1}+V^n}{2}+\beta H^\circ\frac{W^{n+1}+W^n}{2},\tag{3.86}$$

$$U^{n+1}-U^n=2\phi^\circ(\phi^{n+1}-\phi^n),\tag{3.87}$$

$$V^{n+1}-V^n=(\rho^{n+1}-\rho^n)-\boldsymbol{Z}^\circ\cdot\nabla(\phi^{n+1}-\phi^n),\tag{3.88}$$

$$W^{n+1}-W^n=\frac12 H^\circ(\rho^{n+1}-\rho^n),\tag{3.89}$$

*with periodic boundary condition being imposed, where $\phi^\circ=\frac32\phi^n-\frac12\phi^{n-1}$, $\boldsymbol{Z}^\circ=\boldsymbol{Z}(\frac32\phi^n-\frac12\phi^{n-1})$ and $H^\circ=H(\frac32\rho^n-\frac12\rho^{n-1})$.*

| $\delta t$ | LS1 | Order | LS2 | Order |
|---|---|---|---|---|
| $1 \times 10^{-2}$ | 3.97E(-5) | – | 6.03E(-6) | – |
| $5 \times 10^{-3}$ | 1.98E(-5) | 1.00 | 1.51E(-6) | 2.00 |
| $2.5 \times 10^{-3}$ | 9.79E(-6) | 1.02 | 3.80E(-7) | 1.99 |
| $1.25 \times 10^{-3}$ | 4.74E(-6) | 1.05 | 9.48E(-8) | 2.00 |
| $6.25 \times 10^{-4}$ | 2.21E(-6) | 1.10 | 2.33E(-8) | 2.02 |
| $3.125 \times 10^{-4}$ | 1.10E(-6) | 1.00 | 5.53E(-9) | 2.07 |
| $1.5625 \times 10^{-4}$ | 5.25E(-7) | 1.07 | 1.34E(-9) | 2.05 |

TABLE 1. The $L^2$ numerical errors at $t = 0.5$ that are computed by schemes LS1, LS2 using various temporal resolutions with the initial conditions of (4.2), for mesh refinement test in time. The order parameters are of (4.1) and $129^2$ Fourier modes are used to discretize the space.

*The property of symmetric positive definite and unconditional energy stability can be derived thereafter. We omit the details here and leave them to the intersted readers.*

## 4. Numerical experiments

We now present numerical experiments in two dimensions to validate the theoretical results derived in the previous section and demonstrate the efficiency, energy stability and accuracy of the proposed numerical schemes. In all examples, we set the domain $\Omega = [0, 2\pi]^d, d = 2, 3$. If not explicitly specified, the default values of order parameters are given as follows,

$$\widehat{\epsilon} = 0.0001, \quad M_1 = M_2 = 0.01, \quad \epsilon = 0.05, \quad \alpha = 0.01, \quad \beta = 0.05, \quad B = 1. \tag{4.1}$$

We use the Fourier-spectral method to discretize the space, and $129^d$ Fourier modes are used so that the errors from the spatial discretization is negligible compared with the time discretization errors.

4.1. **Accuracy test.** We first test convergence rates of the two numerical schemes, the first order scheme (3.16)-(3.22) (denoted by LS1) and the second order BDF2 scheme (3.57)-(3.63) (denoted by LS2). The following initial conditions

$$\begin{cases} \phi_0(x, y) = 0.1\cos(3x) + 0.4\cos(y), \\ \rho_0(x, y) = 0.2\sin(2x) + 0.5\sin(y) \end{cases} \tag{4.2}$$

are used. We perform the refinement test of the time step size, and choose the approximate solution obtained by using the scheme LS2 with the time step size $\delta t = 7.8125 \times 10^{-5}$ as the benchmark solution (approximately the exact solution) for computing errors. We present the summations of the $L^2$ errors of two phase variables between the numerical solution and the exact solution at $t = 0.5$ with different time step sizes in Table 1. We observe that the schemes, LS1 and LS2, are first order and second order accurate respectively. Moreover, the second order scheme LS2 gives better accuracy than the first order scheme LS1 does when using the same time step.

4.2. **Spinodal decomposition in 2D.** In this example, we study the phase separation behaviors that are called spinodal decomposition using the second order scheme LS2. The process of the phase separation can be studied by considering a homogeneous binary mixture, which is quenched into the unstable part of its miscibility gap. In this case, the spinodal decomposition takes place, which manifests in the spontaneous growth of the concentration fluctuations that leads the system from the homogeneous to the two-phase state. Shortly after the phase separation starts, the domains of the binary components are formed and the interface between the two phases can be specified [2,8,62].

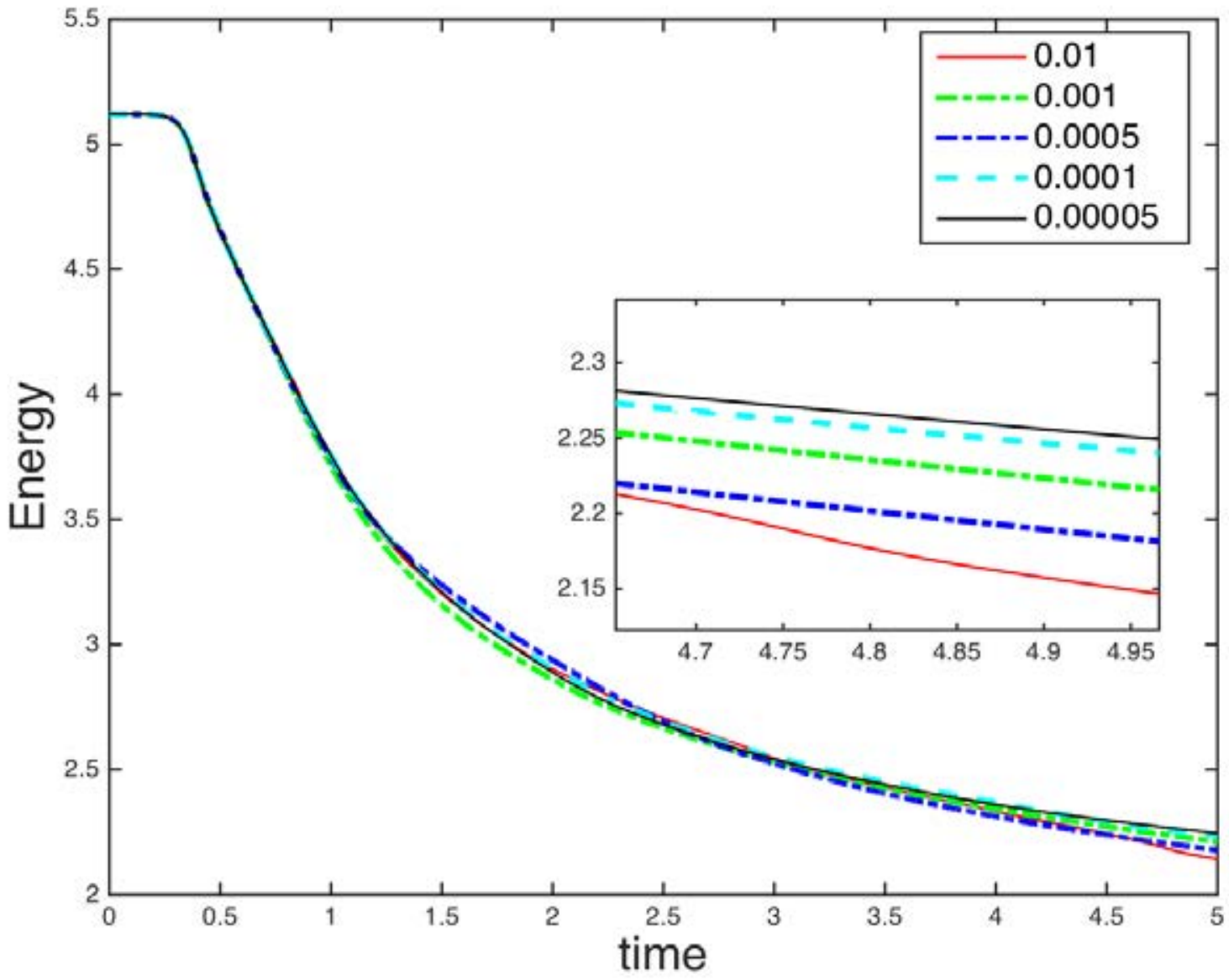

FIGURE 1. Time evolution of the free energy functional for five different time step sizes, $\delta t = 0.00005$, $0.0001$, $0.0005$, $0.001$ and $0.01$ for spinodal decomposition of Example 4.2, with default parameters (4.1). The energy curves show the decays for all time step sizes, that confirms that our algorithm is unconditionally stable.

The initial conditions are taken as the randomly perturbed concentration fields as follows:

$$\phi_0 = \bar{\phi}_0 + 0.001\text{rand}(x, y), \tag{4.3}$$

$$\rho_0 = 0.3 + 0.001\text{rand}(x, y), \tag{4.4}$$

where the rand$(x, y)$ is the random number in $[-1, 1]$ and has zero mean.

We first choose $\bar{\phi}_0 = 0$ and compare the evolution of the free energy functional for five different time step sizes for the default parameters (4.1) until $t = 5$ in Fig. 1. We observe that all five energy curves show the decays for all time step sizes, that confirm that our algorithms are unconditionally stable for any time step. Furthermore, when the time step size $\delta t$ is 0.01, the energy curve is considerable (but not very far) away from others. This means the time step size has to be smaller than 0.01 at least, in order to get reasonably good accuracy. We choose time step $\delta t = 0.0005$ to perform two simulations by varying the value of the initial value $\bar{\phi}_0$ and fixing all other parameters as defined in (4.1).

In Fig. 2, we show the snapshots of coarsening dynamics with $\bar{\phi}_0 = 0$ that means the volume of the immiscible fluids (e.g. oil and water) are the same. Initially, the two fluids are well mixed, and they sooner start to decompose and accumulate. We observe that a relatively high value of the concentration variable $\rho$ is located at the interface. The final equilibrium solution is obtained after $t = 1500$, where the two fluids form a banded shape.

We choose $\bar{\phi}_0 = 0.3$ in Fig. 3. We observe that the fluid component with the less volume accumulates to small satellite drops everywhere. The dynamical behaviors are totally different from last example. When

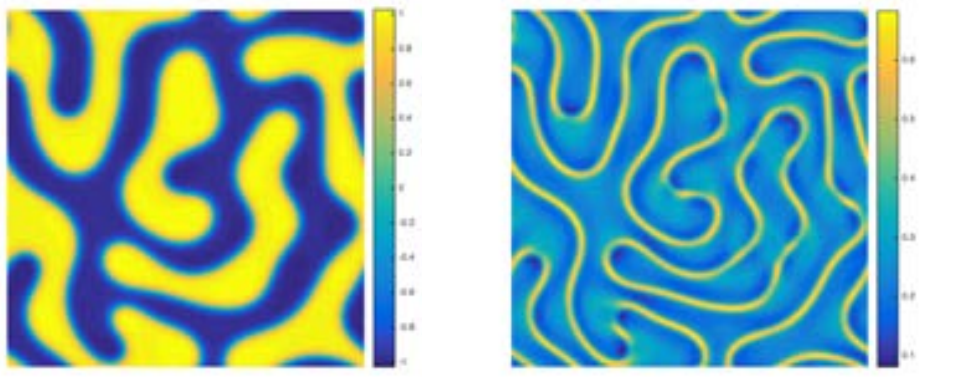

(a) $t = 10$.

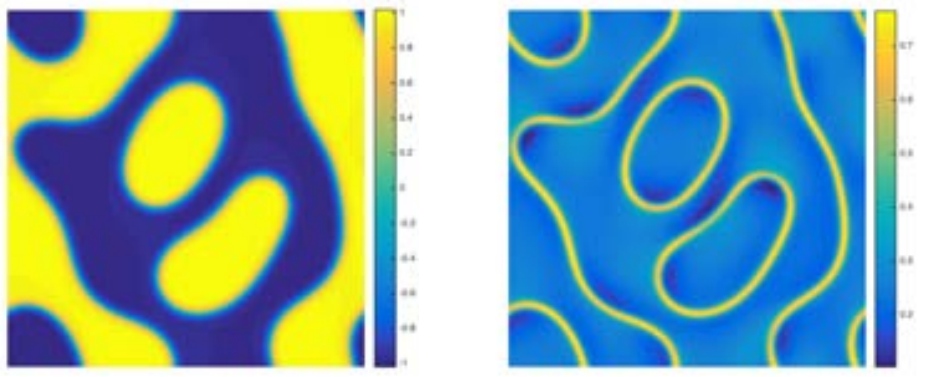

(b) $t = 50$.

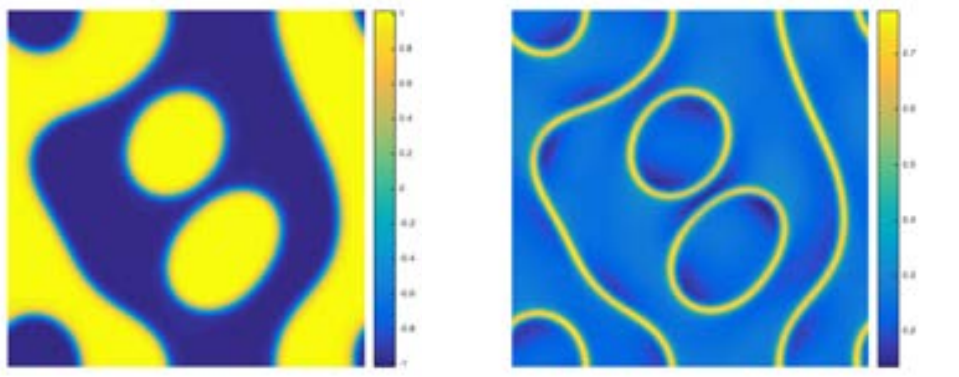

(c) $t = 70$.

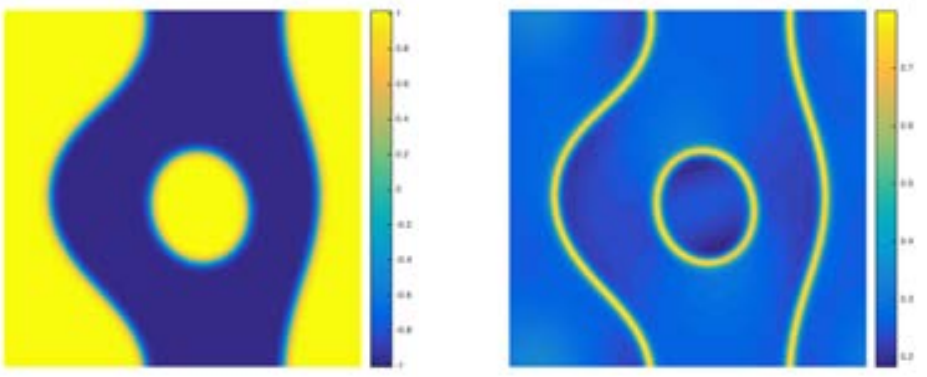

(d) $t = 200$.

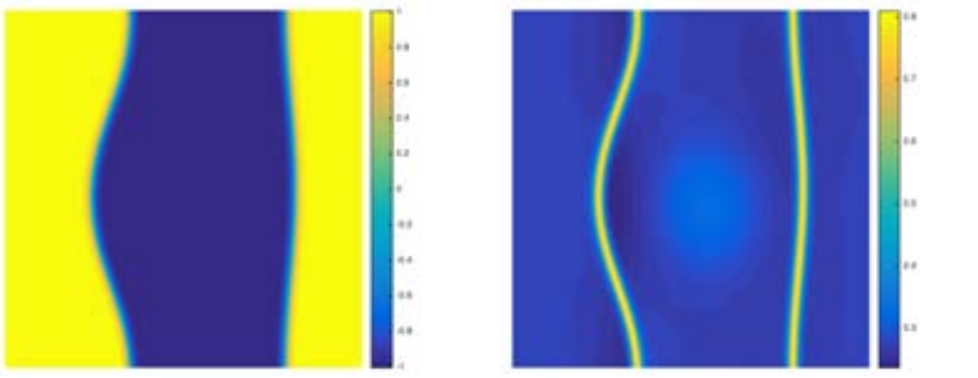

(e) $t = 400$.

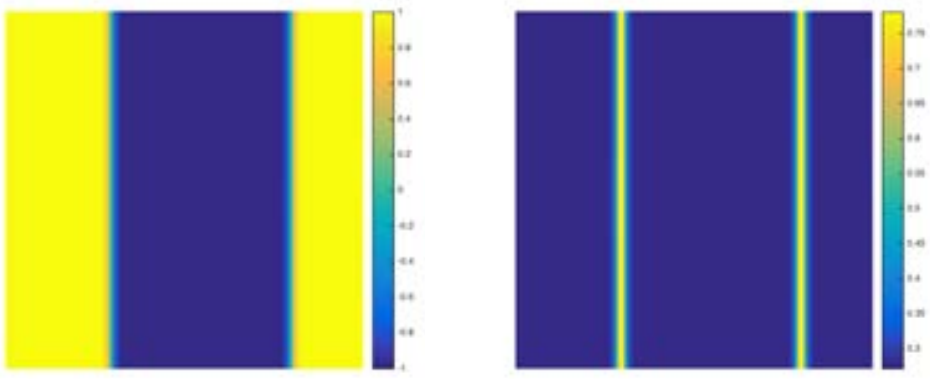

(f) $t = 1500$.

FIGURE 2. 2D spinodal decomposition for random initial data with $\bar{\phi}_0 = 0$. Snapshots of phase variables $\phi$ and $\rho$ are taken at $t = 10$, 50, 70, 200, 400 and 1500. For each subfigure, the left is the profile of $\phi$ and the right is the profile of $\rho$.

the time evolves, the small drops will collide, merge and form drops with bigger size. The final equilibrium solution is obtained after $t = 1500$, where all bubbles accumulate into a big bubble.

In Fig. 4, we plot the evolution of energy curves for both cases that show the energy monotonically decay with respect to the time for both cases.

4.3. **Spinodal decomposition in 3D.** We continue to perform the phase separation dynamics using the second order scheme and time step $\delta t = 0.0005$, but in 3D space. In order to be consistent with the 2D case, the initial condition reads as follows:

$$\phi_0(x, y, z) = \bar{\phi}_0 + 0.001\text{rand}(x, y, z), \tag{4.5}$$

$$\rho_0(x, y, z) = 0.3 + 0.001\text{rand}(x, y, z), \tag{4.6}$$

where the $\text{rand}(x, y, z)$ is the random number in $[-1, 1]$ and has zero mean.

In Fig. 5, we observe the coarsening dynamics with $\bar{\phi}_0 = 0$. As the 2D case, the two fluids initially are well mixed, and they sooner start to decompose and accumulate. But the final steady state (shown in the last subfigure of Fig. 5 is quite different from the 2D result. In Fig. 6, we show the coarsening dynamics for $\bar{\phi}_0 = 0.3$. As the 2D case, the two fluids start to decompose into small drops, and evolve to accumulate.

In Fig. 7, we combine the two steady states with $\bar{\phi}_0 = 0, 0.3$ together to see the difference. In order to obtain more accurate view, since the computed domain is periodic, we then plot the isosurface for 4

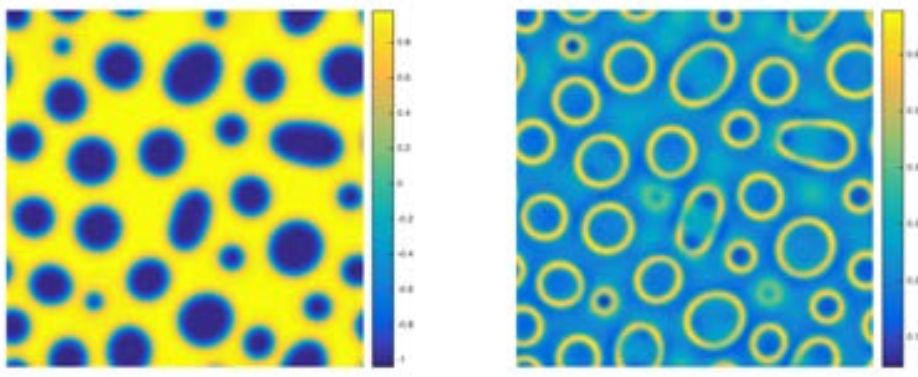

(a) $t = 10$.

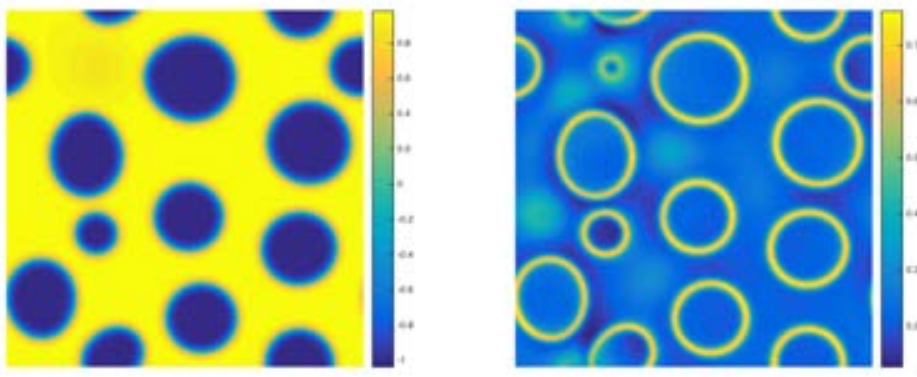

(b) $t = 50$.

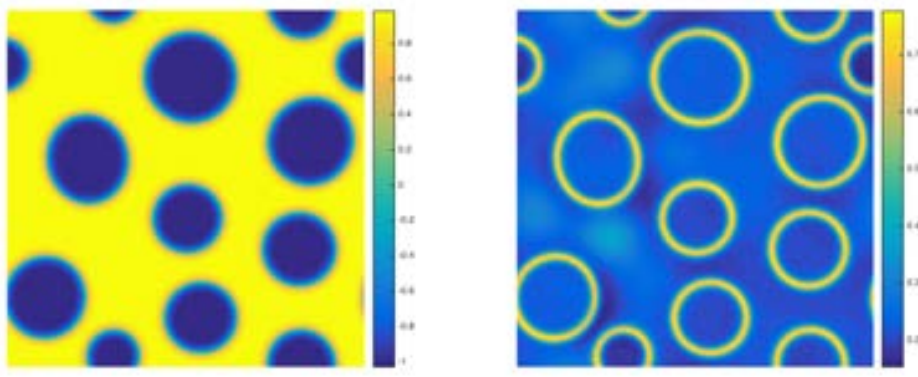

(c) $t = 70$.

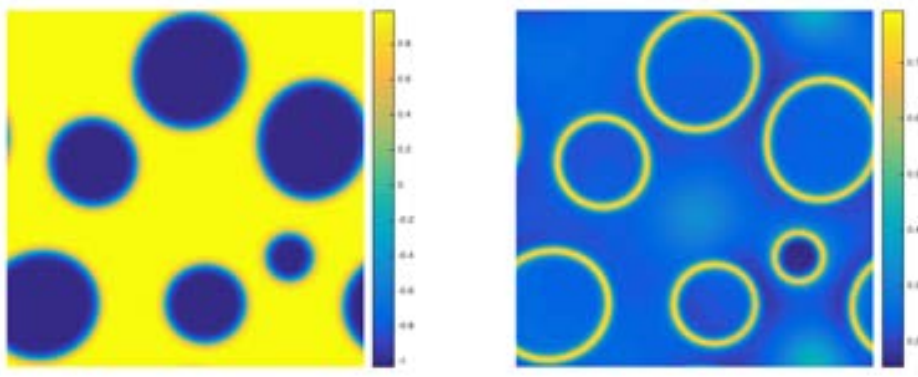

(d) $t = 200$.

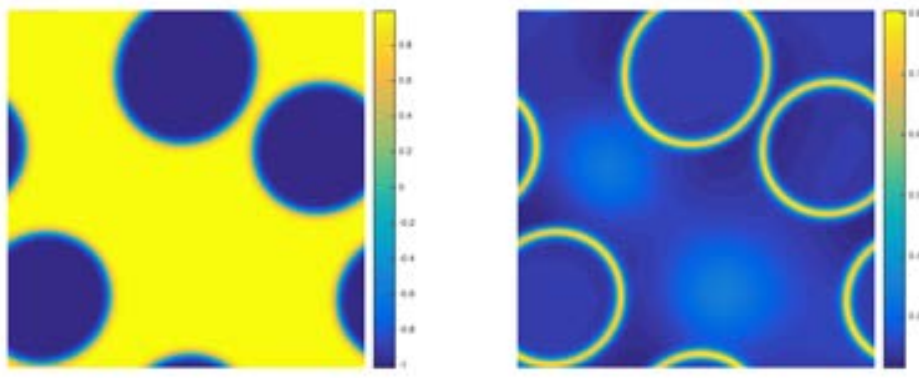

(e) $t = 400$.

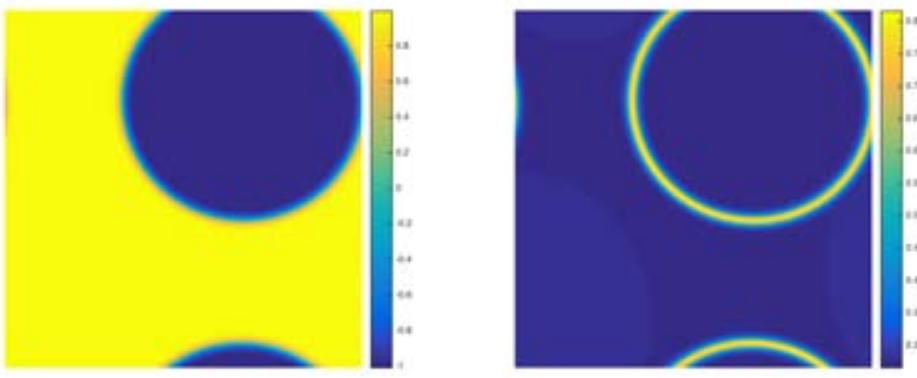

(f) $t = 1500$.

FIGURE 3. 2D spinodal decomposition for random initial data of $\bar{\phi}_0 = 0.3$. Snapshots of phase variables $\phi$ and $\rho$ are taken at $t = 10$, 50, 70, 200, 400, 1500. For each subfigure, the left is the profile of $\phi$ and the right is the profile of $\rho$.

periods, i.e., $[0, 8\pi]$. In Fig. 8, we plot the evolution of energy curves for both cases that show the energy monotonically decay with respect to the time for both cases.

4.4. **Surfactant absorption.**

4.4.1. *Surfactant uniformly distributed initially.* We assume the fluid interface and the surfactant are uniformly distributed over the domain initially and the specific profiles (shown in Fig. 9 (a)) are chosen

$$\phi_0(x, y) = 0.1 + 0.01\cos(6x)\cos(6y), \tag{4.7}$$

$$\rho_0(x, y) = 0.2 + 0.01\cos(6x)\cos(6y). \tag{4.8}$$

We take the time step $\delta t = 0.001$ and use the second order scheme LS2.

Fig. 9 shows the snapshots of coarsening dynamics at $t = 0, 10, 50, 200, 500, 1000$ and the final steady shape forms a large drop. Driven by the coupling entropy energy term, the surfactant is absorbed into the binary fluid interfaces so that the higher concentration appears near the interfaces than other regions.

4.4.2. *Surfactant locally distributed initially.* We here assume the fluid interface and the surfactant field are mismatched over the domain initially. When $t = 0$, the phase field variable $\phi$ is same as the previous example, but the surfactant concentration variable $\phi$ is locally accumulated in the center at the center of the domain.

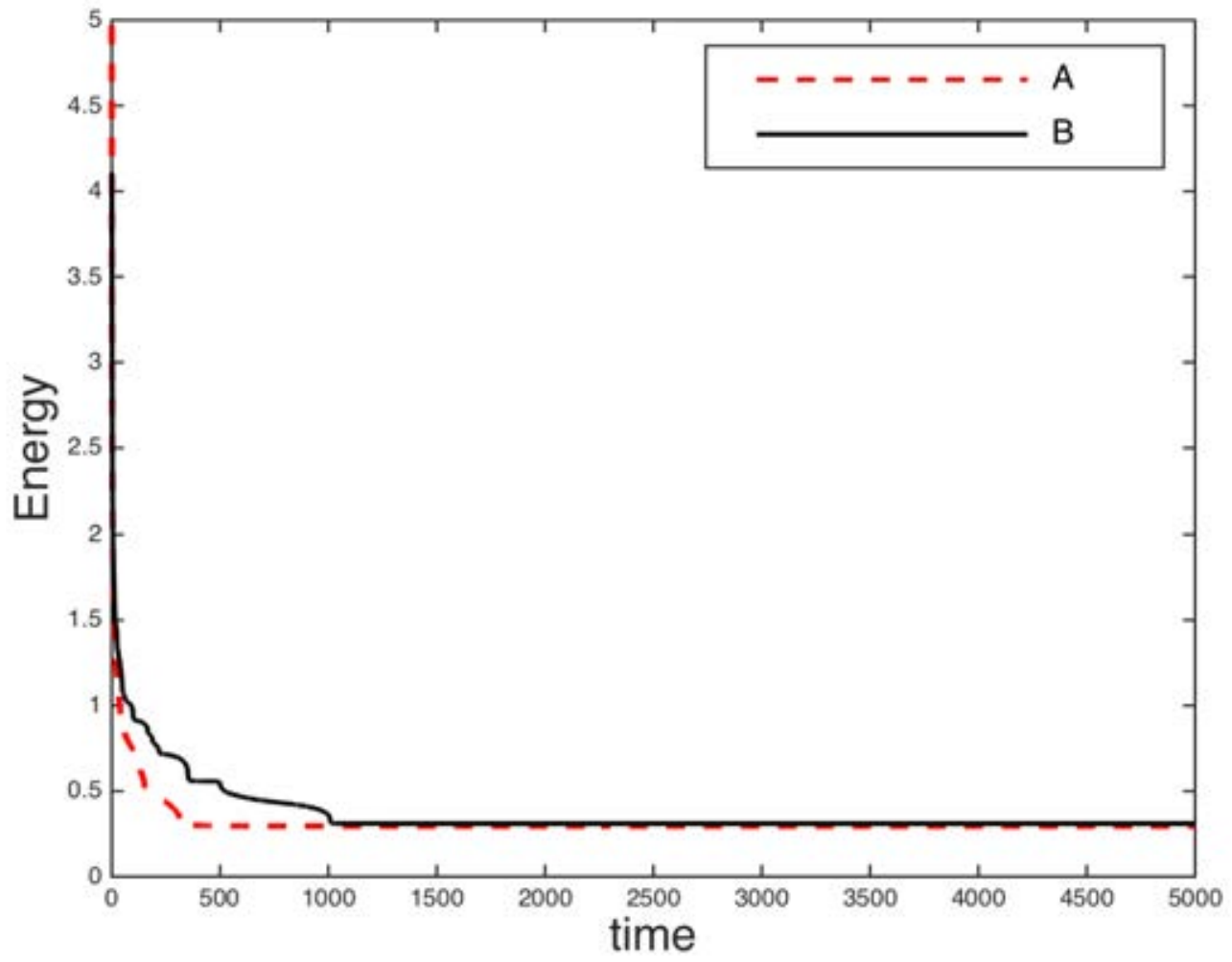

FIGURE 4. Time evolution of the free energy functional for 2D spinodal decomposition Example 4.2 with the initial random data of $\bar{\phi}_0 = 0$ (A) and $\bar{\phi}_0 = 0.3$ (B). The energy curves show monotone decays for both cases.

The initial profiles (shown as Fig. 10 (a)) are chosen as

$$\phi_0(x,y) = 0.1 + 0.01\cos(6x)\cos(6y), \tag{4.9}$$

$$\rho_0(x,y) = 0.8\exp\Big(-\frac{(x-\pi)^2+(y-\pi)^2}{1.25^2}\Big). \tag{4.10}$$

We still take time step $\delta t = 0.001$ and use the scheme LS2 for better accuracy. Fig. 10 shows the snapshots of the dynamics behaviors of $\phi$ and $\rho$ at various times. Since the surfactant is initially concentrated at the center, it takes longer times for the surfactant to diffuse away from this center region. Consequently, during the early stage of the evolution, the higher concentration of surfactant only appears around the center area of the domain. We observe that, the surfactant totally completely diffuses and is absorbed into the binary fluid interfaces at $t = 500$. Overall, the numerical solutions of both examples present similar features to those obtained in [15, 45]. We finally plot the evolution of energy curves in Fig. 11 for Example 4.4.1 and Example 4.4.2, respectively. For both cases, the energy monotonically decays with respect to the time.

## 5. CONCLUDING REMARKS

In this paper, we have developed two efficient, semi-discrete in time, first and second order linear schemes for solving the binary fluid-surfactant phase field model based on a novel IEQ approach. The proposed schemes possess the advantages of the convex splitting method and the stabilized method, and at the same time avoid their defects. More precisely, our schemes are linear, accurate and unconditional energy stable. Furthermore, the induced linear systems are positive definite, thus one could apply any Krylov subspace methods with mass lumping as pre-conditioners to solve them efficiently. We also remark that, to the best of our knowledge, these schemes are not only the first linear and accurate schemes with provable energy stabilities for the particular phase field fluid-surfactant model, but also can serve as a building block to design accurate and stable linear schemes for a class of gradient flow problems with single or multiple variables. Although we consider only time discrete schemes in this paper, the results here can be carried over to any consistent finite-dimensional Galerkin type approximations since the proofs are all based on a variational formulation with all test functions in the same space as the space of the trial functions.

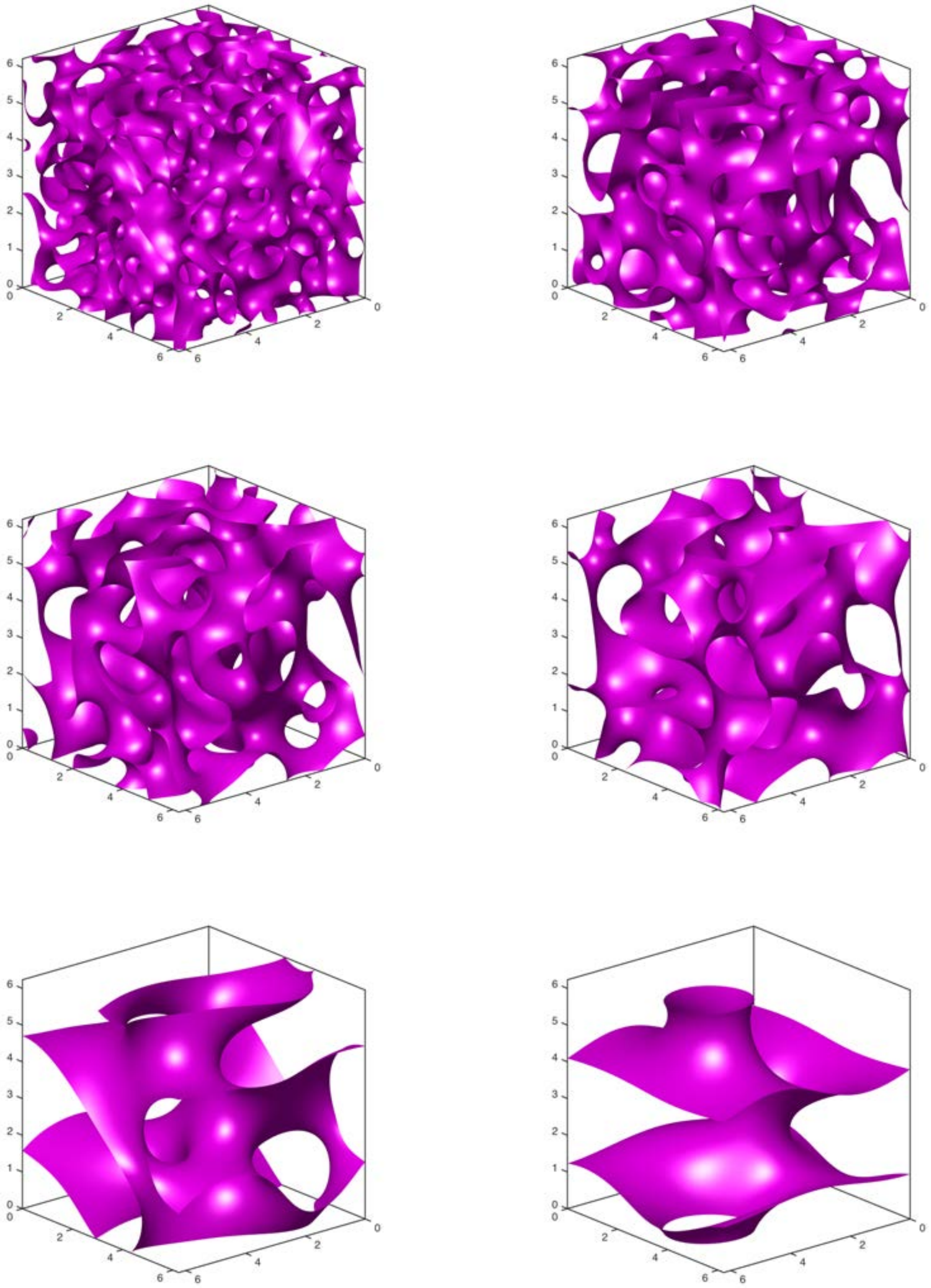

FIGURE 5. 3D spinodal decomposition for $\bar{\phi}_0 = 0$. Snapshots of the phase field variable (isosurface of $\phi = 0$) are taken at $t = 10$, 30, 100, 500, 1000.

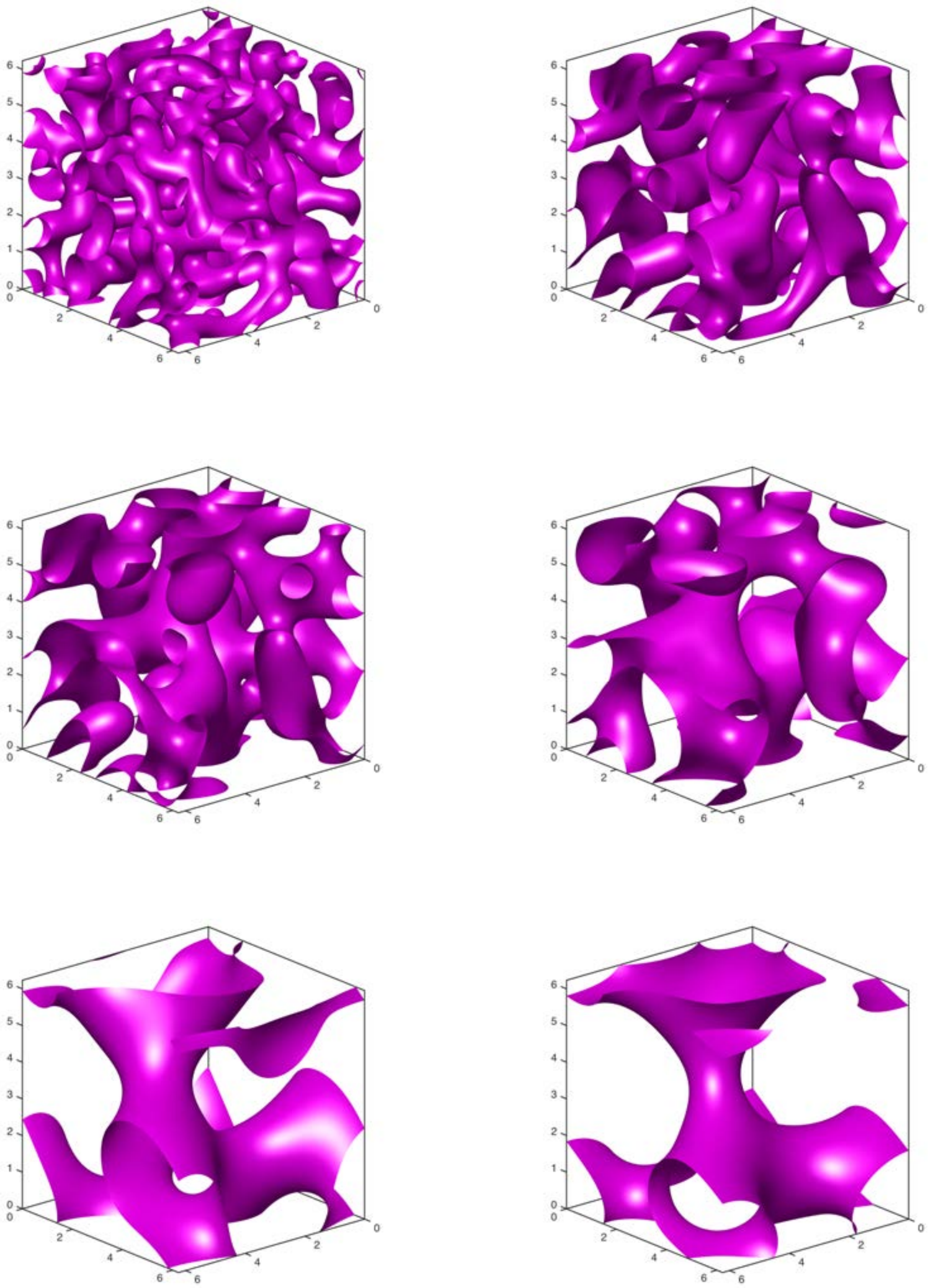

FIGURE 6. 3D spinodal decomposition for $\bar{\phi}_0 = 0.3$. Snapshots of the phase field variable (isosurface of $\phi = 0$) are taken at $t = 10$, 30, 100, 500, 1000.

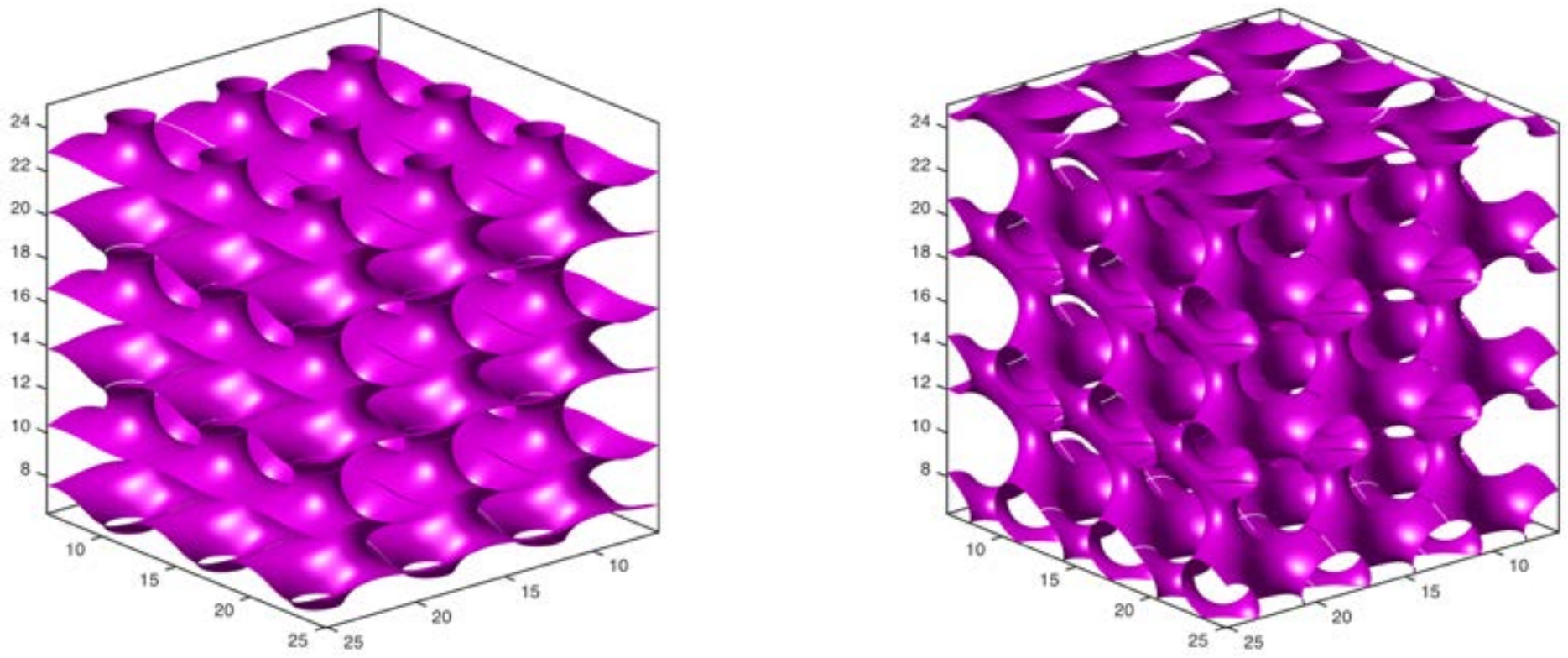

Figure 7. The comparisons of isosurfaces of the equilibrium solutions for 3D spinodal decomposition with two initial values of $\bar{\phi}_0 = 0$ and $\bar{\phi}_0 = 0.3$.

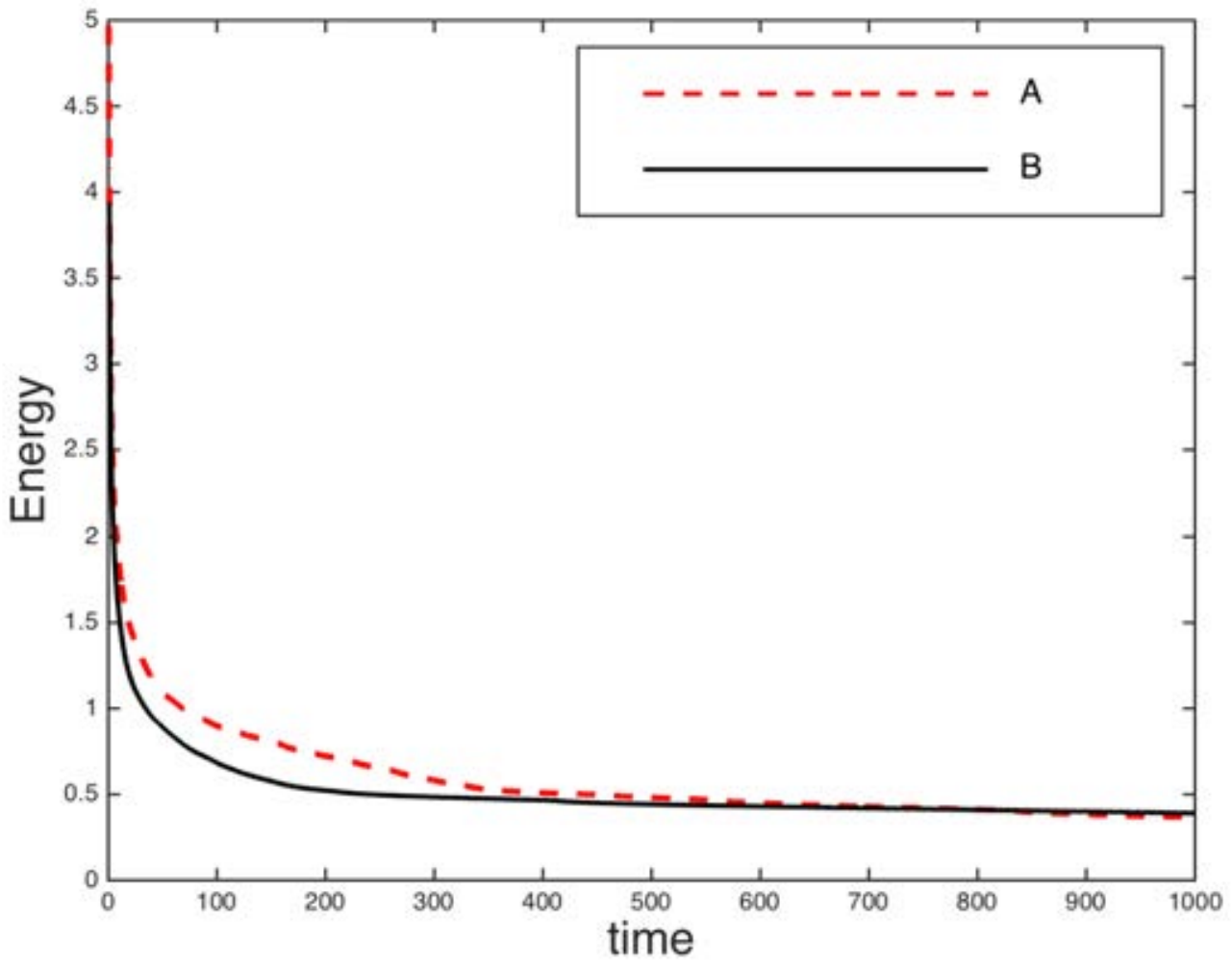

Figure 8. Time evolution of the free energy functional for two three initial conditions of $\bar{\phi}_0 = 0$ (A) and $\bar{\phi}_0 = 0.3$ (B) for the 3D spinodal decompostion. The energy curves show the decays for all time steps, which confirms that our algorithm is unconditionally stable.

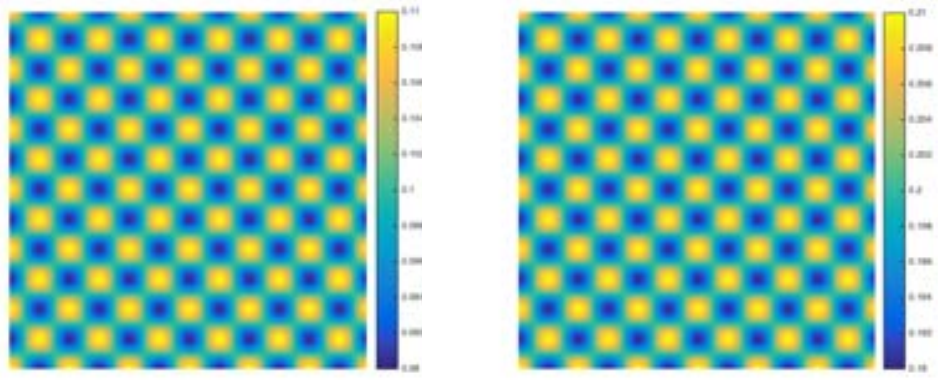

(a) $t = 0$.

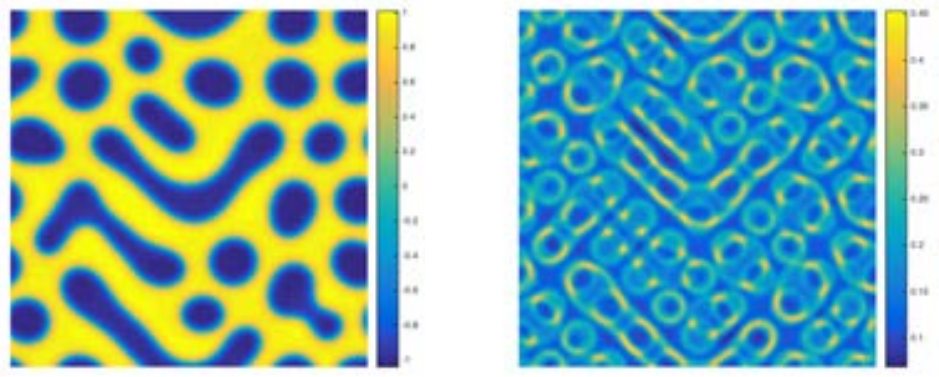

(b) $t = 10$.

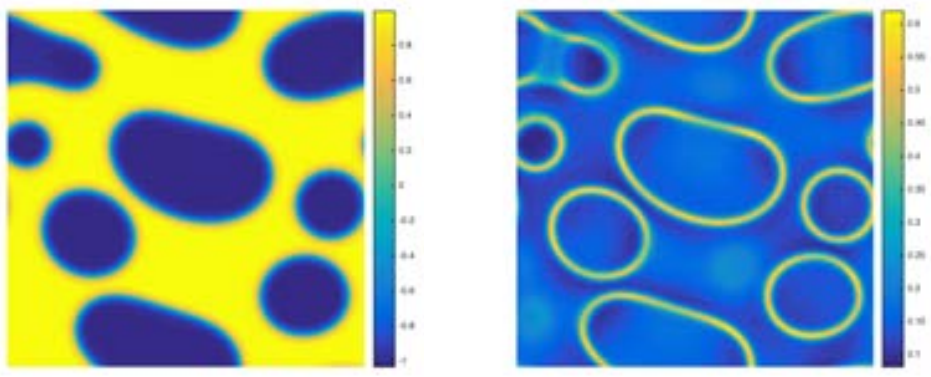

(c) $t = 50$.

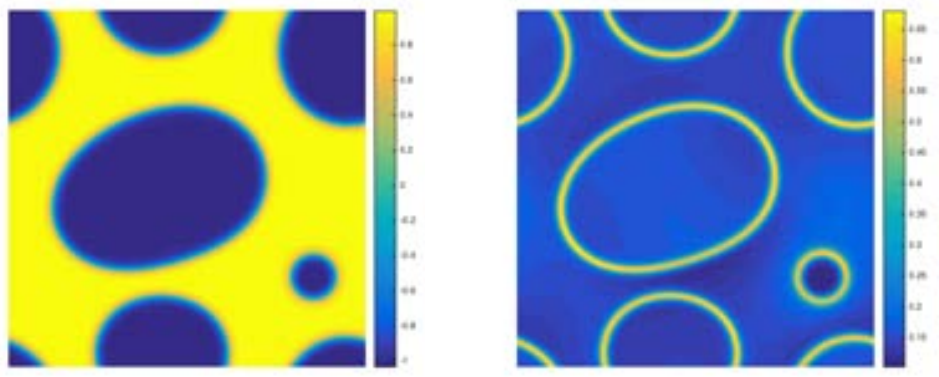

(d) $t = 200$.

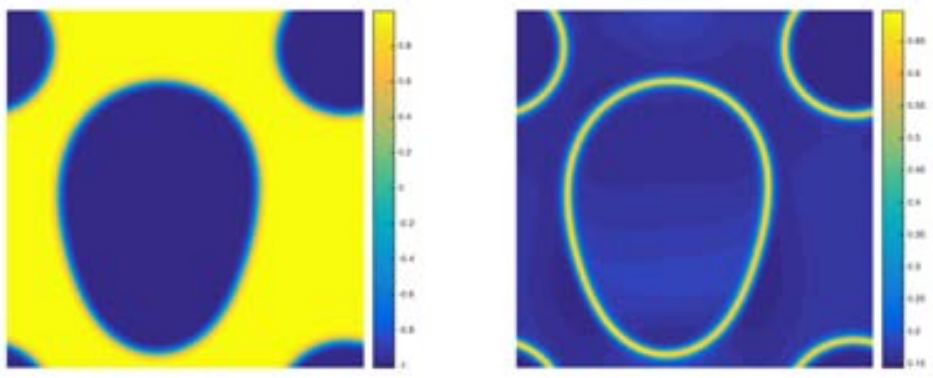

(e) $t = 500$.

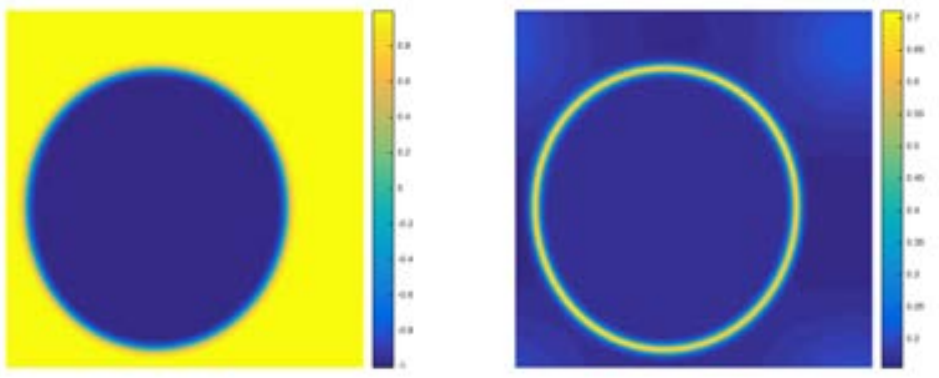

(f) $t = 1000$.

FIGURE 9. Snapshots of the phase variables $\phi$ and $\rho$ are taken at $t = 0$, 10, 50, 200, 500, 1000 for Example 4.4.1, where the surfactants are distributed uniformly initially. For each subfigure, the left is the profile of $\phi$ and the right is the profile of $\rho$.

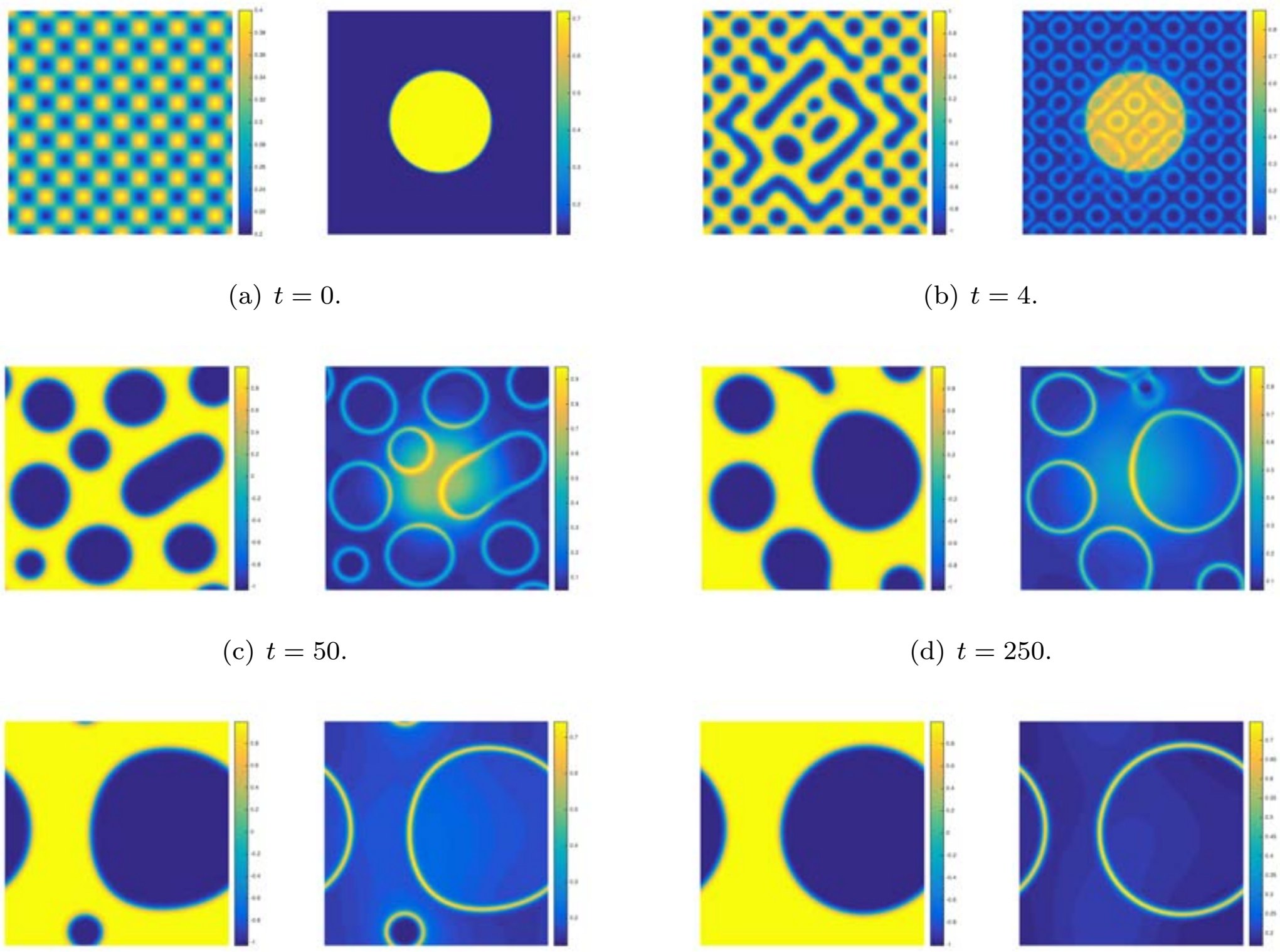

(a) $t = 0$. (b) $t = 4$.

(c) $t = 50$. (d) $t = 250$.

(e) $t = 750$. (f) $t = 1000$.

FIGURE 10. Snapshots of the phase variables $\phi$ and $\rho$ are taken at $t = 0$, 4, 50, 250, 750, 1000 for Example 4.4.2, where the surfactants are distributed uniformly initially. For each subfigure, the left is the profile of $\phi$ and the right is the profile of $\rho$.

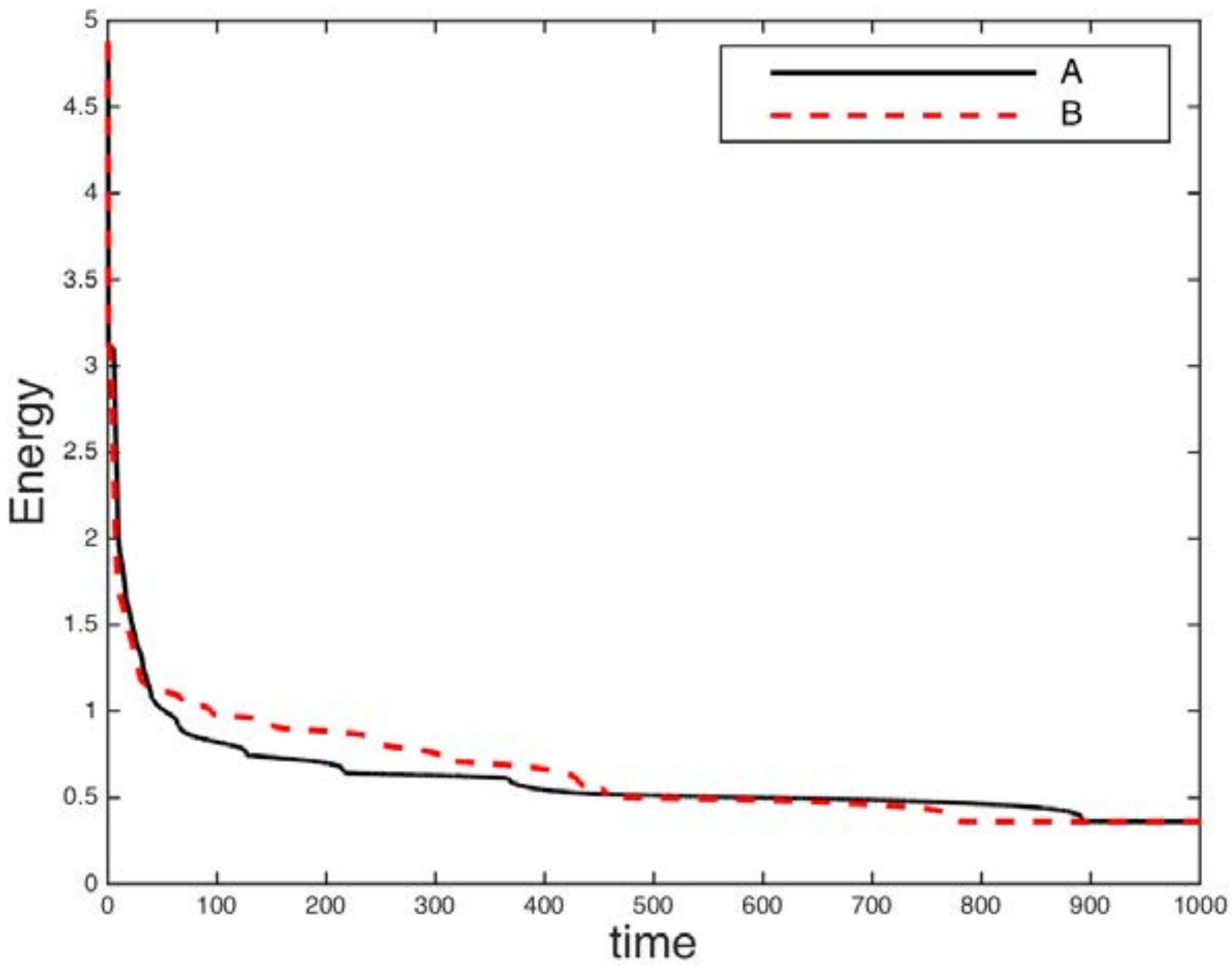

FIGURE 11. Time evolution of the free energy functional for Example 4.4.1 (A) and 4.4.2 (B). The energy curves show monotone decays for both cases.